\documentclass[10pt,a4paper]{article}
\usepackage{stmaryrd} 
\usepackage[pdftex]{graphicx} 
\usepackage{amsmath, amsfonts, amssymb, amsbsy}
\usepackage{latexsym} 

\usepackage[labelfont=bf]{caption}
\usepackage{hyperref}
\hypersetup{
    colorlinks,%
    citecolor=black,%
    filecolor=black,%
    linkcolor=black,%
    urlcolor=black
}

\newcommand{\Div}{{\rm div}}

\newcommand{\ERIC}{\color{black}}

\begin{document}

\newtheorem{theorem}{Theorem}[section]
\newtheorem{proposition}[theorem]{Proposition}
\newtheorem{corollary}[theorem]{Corollary}
\newtheorem{lemma}[theorem]{Lemma}
\newtheorem{definition}[theorem]{Definition}
\newtheorem{remark}[theorem]{Remark}
\newtheorem{example}[theorem]{Example}

\newcommand{\proof}{\noindent \textbf{Proof. }}
\newcommand{\qed}{ \hfill {\vrule width 6pt height 6pt depth 0pt} }

\newcommand{\separe}{\medskip \centerline{\tt -------------------------------------------- } \medskip}

\newcommand{\note}[1]{\medskip \noindent \fbox{\parbox{\textwidth}{\tt #1}} \smallskip}

\newcommand{\eps}{\varepsilon} \renewcommand{\det}{\mathrm{det}} \newcommand{\argmin}{ \mathrm{argmin} \,}
\newcommand{\Om}{\Omega} \def\interior{\mathaccent'27} 
\newcommand{\weakto}{ \rightharpoonup}  \newcommand{\weakstarto}{\stackrel{*}{\rightharpoonup}}
\newcommand{\R}{\mathbb{R}}
\newcommand{\U}{\mathcal{U}}\newcommand{\F}{\mathcal{F}}\newcommand{\E}{\mathcal{E}}\newcommand{\K}{\mathcal{K}}

\newcommand{\stress}{\sigma} \newcommand{\strain}{\epsilon} 
\newcommand{\neu}{ \partial_{\mbox{\it \tiny N}} } \newcommand{\dir}{\partial_{\mbox{\it \tiny D}} }
\newcommand{\I}{ {\mbox{\tiny \rm I}} } \newcommand{\II}{ {\mbox{\tiny \rm II}} }
\newcommand{\jumpo}[1]{\llbracket #1 \hspace{0.5pt}\rrbracket}
\def\Xint#1{\mathchoice
 {\XXint\displaystyle\textstyle{#1}}%
 {\XXint\textstyle\scriptstyle{#1}}%
 {\XXint\scriptstyle\scriptscriptstyle{#1}}%
 {\XXint\scriptscriptstyle\scriptscriptstyle{#1}}%
 \!\int}
\def\XXint#1#2#3{{\setbox0=\hbox{$#1{#2#3}{\int}$}
 \vcenter{\hbox{$#2#3$}}\kern-.5\wd0}}
 \def\ddashint{\Xint=} 
 \def\dashint{\Xint-}

\pagebreak 

\thispagestyle{empty} 


\vspace{0.5cm}
\noindent {\Large \bf Existence, energy identity and higher time regularity of solutions}%
\vspace{6pt}

\noindent{\Large \bf to a dynamic visco-elastic cohesive interface model} 

\vspace{36pt}

\begin{small}
{\bf M.~Negri}

{Department of Mathematics -  University of Pavia} 

{Via A.~Ferrata 1 - 27100 Pavia - Italy}

{matteo.negri@unipv.it} 

\medskip
{and} 
\medskip

{Institute for Applied Mathematics and Information Technologies - CNR}

{Via A.~Ferrata 5 - 27100 Pavia - Italy}

\vspace{24pt}
{\bf R.~Scala}

{Dipartimento di Ingegneria dell'Informazione e Scienze Matematiche - University of Siena}

{San Niccol\`o, via Roma, 56 - 
53100 Siena - Italy}

{riccardo.scala@unisi.it}

\vspace{36pt}
\noindent {\bf Abstract.}  We study the dynamics of visco-elastic materials coupled by a common cohesive interface (or, equivalently, {two single domains separated by} a prescribed cohesive crack) in the anti-plane setting. We consider a general class of traction-separation laws featuring an activation threshold on the normal stress, softening and elastic unloading. In strong form, the evolution is described by a system of PDEs coupling momentum balance (in the bulk) with transmission and Karush-Kuhn-Tucker conditions (on the interface). We provide a detailed analysis of the system.  We first prove existence of a weak solution, employing a time discrete approach and a regularization of the initial data. Then, we prove our main results: the energy identity and the existence of { solutions} with acceleration in $L^\infty (0,T; L^2)$.



\end{small}

\section{Introduction}

The problem of interaction between two elastic or r visco-elastic bodies in contact along a common interface arises naturally in many branches of the mechanics of solid and is strongly connected with the nucleation and propagation of cracks along a prescribed path \cite{Barenblat}. 
In the last few years, these problems have been studied extensively from the mathematical point of view, considering several traction-separation laws and different evolution, both in quasi-statics (energetic and balanced-viscosity) and dynamic setting.

As far as quasi-static evolutions several contributions have been proposed, covering many similar models with either adhesive or cohesive contact. 
For energetic evolutions \cite{MielkeRoubicek} we quote \cite{DMZ,RSZ,CagnettiToader,RMP}, among the many in our context.
%
For balanced-viscosity evolutions we may further distinguish between solutions obtained by vanishing viscosity \cite{MielkeRossiSavare09} and solutions obtained by incremental local minimization \cite{Negri14}. For our problem, we mention \cite{Cagnetti08,Almi} for the latter approach, and \cite{NegriVitali_IFB18, NScala_NARWA17} for the former.


Besides rate-independent evolutions also rate-dependent systems with inertial effects have been considered to describe the motion of visco-elastic bodies together with crack propagation, delamination, debonding and damage evolution; we quote only some important contributions, as \cite{DalMasoLarsen_AANLCSFMNRLMA11,LS2014,DMLN2016,DMLT2016,ScalaSchimperna_EJAM17,Scala_IFB17,DML2017,TZ,RN2020} and references therein. In some cases, starting from dynamic models it is possible to recover a quasi-static evolution by time rescaling. 
Among the results in this direction, we quote \cite{Roubicek_SJMA13,Scala_ESAIMCOCV17} for delamination models, \cite{LN2019} for a debonding problem, \cite{DMS2014} for a plasticity, \cite{LRTTN2018} for damage,  and finally \cite{Ago2012,Nar2017} for general results in finite dimension.

\medskip

Let us now describe the model. Following the quasi-static models of  \cite{OrtizPandolfi_IJNME99,NegriVitali_IFB18} and the dynamical model of \cite{TZ} we are here interested in an evolution accounting for visco-elastic and inertial effects, in the bulk, and traction-separation laws with finite activation threshold and different loading-unloading regimes, in the cohesive interface. A closer comparison with the setting and the results of \cite{TZ} is postponed at the end of the introduction.
 Let $\Om\subset\R^2$ represent (the planar section of) our reference configuration, consisting of a couple of elastic bodies in contact along a common interface $K$, i.e., $\Om=\Om^+\cup\Om^-$ and $K=\partial\Om^+\cap\partial\Om^-$. We denote by $u:\Om\rightarrow\R$ the antiplane displacement, i.e., orthogonal to the plane containing $\Om$. The elastic energy of the system is then written as
\begin{align}
 \mathcal{E} (u) = \int_\Om\mu|\nabla u|^2dx,
\end{align}
where $\mu>0$ is the shear modulus. We consider external loadings in the form of a linear functional $f(t)$, depending on time $t \in [0,T]$, accounting for bulk forces and traction on the Neumann part of the boundary, i.e.,
$$
	( f(t) , u ) = \int_\Omega f_b u \, dx + \int_{\partial_N \Omega} f_s u \, d \mathcal{H}^1.
$$
In addition we take into account a damping term given by the Kelvin-Voigt viscosity, associated to the dissipation
\begin{align}
 \mathcal{D} ( \dot u) = \int_\Om\eta|\nabla \dot u|^2dx .
\end{align}
%
The energy is complemented with a cohesive potential of the form 
\begin{align}
 \Psi(u,\xi):=\int_K\psi(\jumpo{u},\xi)d\mathcal H^1. 
\end{align}
The density $\psi$ depends on the opening $\jumpo{u}:=u^+-u^-$ (the difference between the traces of $u$ on $K$ from $\Om^+$ and $\Om^-$) and on the internal variable $\xi$, whose role is to keep track of the history of the evolution of $\jumpo{u}$: the (non-negative) value of $\xi$ at a point $l \in K$ and at time $t\in[0,T]$ plays the role of the maximum opening at the point $l$ during the time interval $[0,t]$, i.e.
$
 \xi(l,t):= \max_{s\in[0,t]}\{|\jumpo{u(l,s)}|\}.
$
Besides some technical assumptions, the density $\psi$ is monotone non-decreasing in the internal variable and convex-concave in the opening.  
The former property models dissipation on the interface while the latter models the elastic response in unloading and softening in loading. In particular, for $\xi=0$ softening occurs only above a certain threshold on the normal traction, indeed the density $\psi$ behaves like $ c |\jumpo{u}|$ for $\xi=0$ and $| \jumpo{u}| \ll 1$ and thus the energy $\Psi$ turns out not to be Frech\'et differentiable.
%

In summary, the total energy $\mathcal F$ of our system is given by
\begin{align}
\mathcal F(t,u,\xi): = \mathcal E(u) + \Psi(u,\xi) - ( f(t) , u ) .
\end{align}
Along with these energies, we study the dynamics of the following system, consisting of the momentum balance
\begin{subequations}
\begin{align}\label{equazione1}
\rho \ddot u-\Div \sigma=f_b \text{ in }\Om,
\end{align}
where $\sigma:=\mu\nabla u+\eta\nabla \dot u$ is the visco-elastic stress, the transmission condition
\begin{align}\label{equazione2}
 \sigma^+\nu=\sigma^-\nu  \in  \partial_w\psi(\jumpo{u},\xi) \quad\text{ on }K,
\end{align}
where $\nu$ is the normal unit vector to $K$ pointing from $\Om^+$ to $\Om^-$, $\partial_w\psi$ represents the subdifferential of $\psi$ with respect to $w=\jumpo{u}$, and the flow rule for the internal variable
\begin{align}\label{equazione3}
 \dot \xi(\xi-|\jumpo{u}|)=0,\;\;|\jumpo{u}|\leq \xi\quad\text{ on }K.
\end{align}
Equations \eqref{equazione1}-\eqref{equazione3} are complemented with the boundary conditions 
\begin{align}\label{equazione5}
 &u=0\text{ on }\partial_D\Om, \\ &  \sigma\nu = f_s \text{ on }\partial_N\Om,
\end{align}
and the initial condition 
\begin{align}\label{equazione6}
 u(0)=u_0,\;\;\;\;\dot u(0)=v_0, \;\;\;\;\xi(0)=\xi_0\geq|\jumpo{u_0}|.
\end{align}
\end{subequations}
The system \eqref{equazione1}-\eqref{equazione6} provides a strong formulation; {existence of solutions actually requires some additional conditions on the energies and on the data.}
{Instead, under more general hypotheses we prove the energy identity and the existence of weak solutions, in the following sense.}
%
Assuming $u_0,v_0\in H^1(\Om)\times L^2(\Om)$, we prove that there exist  $u$, $\xi$ with 
\begin{align}
 &u\in H^1(0,T;H^1(\Om))\cap H^2(0,T;H^1(\Om)^*),\nonumber\\
 &\xi\in H^1(0,T;L^2(K)),
\end{align}
satisfying \eqref{equazione5}, \eqref{equazione6}, and for a.e. $t\in[0,T]$
\begin{align}\label{eq_debole1}
 &( \rho\ddot u(t),\phi)_{H^1} +\partial_u\mathcal F(t,u(t),\xi(t);\phi)+\langle\eta\nabla \dot u,\nabla \phi\rangle \geq0,\\
 &\dot \xi(t)(\xi(t)-|\jumpo{u(t)}|)=0,\;\;|\jumpo{u(t)}|\leq \xi(t)\quad \text{ on $K$},\label{eq_debole2}
\end{align}
for every $\phi\in H^1(\Om)$ with $\phi=0$ on $\partial_D\Om$. In equation \eqref{eq_debole1} the symbol $\partial_u\mathcal F(t,u,\xi;\phi)$ represents the partial derivative of the total energy in the direction $\phi$, whereas $\langle\cdot,\cdot\rangle$ is the scalar product in  $L^2(\Om)$, and $( \cdot ,\cdot )_{H^1}$ is the duality between $H^1(\Omega)^*$ and $H^1(\Omega)$.
We prove that the solutions to \eqref{eq_debole1} satisfy,  for any $0\leq s\leq t\leq T$, the energy balance
\begin{align}\label{bilancio}
	\E ( u( t) ) + \Psi( u(t) , \xi (t)) + \K ( \dot{u} (t) )  & =  \, \E ( u (s)) + \Psi ( u(s), \xi(s)) + \K ( v(s) ) \nonumber \\ & \quad + \int_{s}^{t} (f( r) , \dot{u} (r) ) \, dr - \int_{s}^{t} \mathcal{D} ( \dot{u} (r) )  \,dr ,  
\end{align}
%
where 
$$
  \mathcal{K}( \dot{u} ) = \int_\Omega | \dot{u}|^2 \, dx 
$$
is the kinetic energy. We also show that these solutions solve the system of equations \eqref{equazione1}-\eqref{equazione6} if some additional technical requirements are made. {More precisely, if the domain, the initial data, and the external force, are suitably regular, we show that \eqref{equazione1}-\eqref{equazione6} hold with $\ddot u$ belonging to $L^2(\Omega)$. These are the so-called strong solutions.}

The precise statements of our main results, namely Theorems \ref{t.teor} and \ref{d.wsol2}, are contained in Section \ref{mainresults_sec}. In order to prove them, we proceed in two steps. The idea consists in fixing a parameter $\eps>0$ and finding a solution $(u^\eps,\xi^\eps)$ which satisfies \eqref{eq_debole1}-\eqref{eq_debole2} with regularized initial data. In particular the initial values of the internal variable $\xi^\eps (0)$ is strictly positive. In this way we gain a regularized cohesive energy (since the cohesive potential is singular only at $\xi=0$). In a second step we pass to the limit as $\eps\rightarrow0$ and recover the original initial data and the solution to the original problem. In order to prove the existence of $(u^\eps,\xi^\eps)$ we classically proceed by time discretization and solve a minimization problem at each step. 

{We conclude this introduction comparing with adhesive interface energies (more details are provided in Remark \ref{comparison}); as pointed out in \cite{TZ}, adhesive and cohesive settings are indeed closely related,  with some differences which are worth to point out. 
Although adhesive models cover a vast type of interface energy profiles, see e.g.~\cite{TZ}, the interface energy is, roughly speaking, quadratic at the origin while the dissipation on the interface is  positively $1$-homogeneous. Here we employ instead more general cohesive energies, non-quadratic and non-differentiable in the origin, combined with non-linear dissipations; as a consequence our model features a finite positive activation threshold on the normal tractions, typical of cohesive fracture. This provides additional mathematical difficulties which enforces us to first regularize the problem  and then to pass to the limit  to recover a solution to the original singular model. The regularization, performed indirectly with the aid of the initial values $\xi^\eps (0)$, allows in practice to work with a family of adhesive energies, becoming singular as $\eps \to 0$. 

It is also interesting to briefly compare with the so-called \textit{semi-stable energetic solutions}, often employed in adhesive models. This notion of solution, in essence, consists in coupling \eqref{eq_debole1} with an energy minimization problem, governing the evolution of the internal variable (see also Remark \ref{comparison}). In our context, recasting the discrete scheme as energy minimization problem for the internal variable leads to an ill-posed problems, with infinitely many solutions. For this reason, we employ instead the Karush-Kuhn-Tucker conditions \eqref{eq_debole2} which describe clearly the evolution of the internal variable and are crucial in the proof of the energy identity.
}

Beyond existence of weak solutions, our work aims at proving a couple of noteworthy properties, which have never been proved in the literature on cohesive interfaces: the energy identity \eqref{bilancio} and the existence of strong solutions, i.e.~solutions to the PDE system \eqref{equazione1}-\eqref{equazione6}. Both the properties are natural from the modelling point of view, however their proofs are not straightforward. Energy identity requires indeed to show that the interface energy is absolutely continuous in time, which in turn follows from the time regularity of the internal variable together with the Karush-Kuhn-Tucker conditions; strong solutions require to improve the time regularity of the acceleration $\ddot{u}$  from $L^2(0,T; H^1(\Omega)^*)$ to $L^\infty (0,T; L^2(\Omega))$ by means of a delicate compactness estimate, which, loosely speaking, holds when the energy  $u \mapsto \F ( t , u , \xi)$ is convex and boundary traction vanishes.

\pagebreak

\tableofcontents

\section{Setting and preliminaries\label{s.setting}}

\subsection{Geometry and spaces}



We consider as reference configuration a set $\Omega \subset \R^2$ obtained as union of disjoint open sets $\Omega^+$ and $\Omega^-$, whose boundary overlaps on a common interface $K = \partial \Omega^+ \cap \partial \Omega^-$. Technically, we assume both the sets $\Omega^\pm$ to be connected, bounded, and Lipschitz (additional  hypotheses on $\Omega$ and $K$ will be introduced to prove existence of strong solutions).
We assume that the sets $\partial \Om^\pm \setminus K$ (i.e., the set where the boundaries $\partial \Omega^\pm$ are not in contact) can be written as the union of two parts, $\partial_D \Om^\pm$ and $\partial_N \Om ^\pm = \partial \Omega^\pm \setminus (K \cup \partial_D \Omega^\pm)$.
On the former part of the boundary we will impose Dirichlet boundary condition, whereas in the latter Neumann boundary condition. In order to avoid trivial and singular cases, which would require ad-hoc arguments, we will assume that $\mathcal{H}^1 (K) >0$ and $\mathcal{H}^1(\partial_D \Omega^\pm)>0$. 
For later convenience we will also denote $\partial_D \Om = \partial_D \Om^+ \cup \partial_D \Om^-$ and similarly $\partial_N \Om =  \partial_N \Om^+ \cup \partial_N \Om^-$. Note that 
$$  \partial \Om = K \cup \partial_D \Om \cup \partial_N \Om .$$

We denote by $\nu^\pm$ the outer unit normal to $\Om^{\pm}$; we also denote by $\nu$ the outer unit normal on $\partial \Omega\setminus K$. On $K$ we set $\nu:= \nu^-=-\nu^+$. 



We consider antiplane displacements $u : \Om \to \R$ belonging to the space
$$\U = \{ u \in H^1(\Omega) : u =0 \text{ on $\partial_D \Omega$} \}. $$
The space $\U$ is endowed with the norm $\| u \|_{\U} =\| \nabla u \|_{L^2} $; by Poincar\'e inequality in the sets $\Omega^\pm$ this norm is equivalent to $\| u \|_{H^1} = \|u\|_{L^2} +\| \nabla u \|_{L^2}  $. 
We denote by  $\U^*$ the dual space of $\U$, and by $( \cdot , \cdot )_\U$ the duality pairing between $\U^*$ and $\U$; we denote the $L^2$-scalar product by $\langle \cdot, \cdot \rangle$.


We denote by $u^\pm$ the restriction to $K$ of the traces $u^\pm$ on the Lipschitz boundaries $\partial \Om^\pm$. The jump of $u$ on $K$ is $\jumpo{u} = u^+ - u^-$. Clearly the map $u \mapsto \jumpo{u}$ is linear and continuous from $\U$ to $L^2(K)$. 
By Poincar\'e and trace theorem we have the following inequality
\begin{equation} \label{ineq_trace}
     \| \nabla u \|^2_{L^2(\Omega)} \ge \hat c \int_K \jumpo{ u}^2 \, d\mathcal{H}^1, 
\end{equation}
for some positive constant $\hat c$.
On the interface we will also employ an internal (history) variable $\xi$, which will play the role of the maximal opening, belonging to the space $$\Xi = L^2(K; \R_+),$$ where $\R_+ = [0, +\infty)$; the set $\Xi$ is endowed with the $L^2$-norm.

\subsection{Energies}



First, let us introduce the energy $\F : [0,T] \times \U \times \Xi \to \mathbb{R}$ of the form $$\F ( t ,u , \xi) = \E ( u )  - (f(t) , u)_{\U} + \Psi ( u , \xi)   , $$
where $\E$ is the elastic (bulk) energy, $f$ accounts for bulk and surface external forces while $\Psi $ is the cohesive interface energy, defined in the next subsection (note that $\Psi$ will further split into the stored energy $\Psi_s$ and the dissipated energy $\Psi_d$). Moreover, we will employ a kinetic energy $\mathcal{K}: \U \to \R$ and a dissipation (rate of dissipated energy) $\mathcal{D}: \U \to \R$ of the form 
$$
   \mathcal{K} (\dot{u}) , \qquad \mathcal{D} (\dot{u} ) = \partial_{v} \mathcal{R} ( \dot{u} ) [ \dot{u} ] ,
$$
where $\mathcal{R} : \U \to \R$ is a dissipation pseudo-potential.

More precisely, the elastic energy reads 
$$
 \E ( u ) = \tfrac12 \int_{\Omega} \mu  | \nabla u |^2 \, dx , 
$$
where $\mu = \mu^+ 1_{\Om^+} + \mu^- 1_{\Om^-}$ and $\mu^\pm >0$ are the shear moduli of $\Omega^\pm$.

In the sequel we will assume $f \in W^{1,2}(0,T; \U^*)$. This hypothesis covers the case 
\begin{align}\label{dec_f}
(f(t),u)_\U=\int_{\Om} f_b (t) \, u\;dx+\int_{\partial_N\Om} f_s (t)\,u \, d \mathcal{H}^1,
\end{align}
where $f_b \in W^{1,2}(0,T;L^2(\Om))$ and $f_s \in W^{1,2} (0,T;L^2(\partial_N\Om))$ are the external forces, acting respectively in the bulk and in the Neumann part of the boundary; 
in order to prove that $\ddot u$ belongs to $L^2(\Om)$ we will actually assume that $f_s=0$ and thus $f\in W^{1,2}(0,T;L^2(\Om))$.
Accordingly, we introduce the power of external forces $\mathcal{P}_{ext} : [0,T] \times \U \to \R$ given by 
$$
	\mathcal{P}_{ext} ( t, \dot{u} ) = ( f(t) , \dot{u} )_\U  .
$$
%
%

Next, given the density $\rho = \rho^+ 1_{\Om^+} + \rho^- 1_{\Om^-}$, for $\rho^\pm >0$, the  kinetic energy $\mathcal{K} : \U \to \R$ is defined by
$$ \mathcal{K} ( v ) = \tfrac12 \int_\Omega \rho | v |^2 dx , $$
where $v$ plays the role of the speed $\dot{u}$.

Then, for $\eta = \eta^+ 1_{\Om^+} + \eta^- 1_{\Om^-}$, with $\eta^\pm >0$, we consider the Kelvin-Voigt visco-elastic dissipation $\mathcal{D} : \U \to \R$ and the associated pseudo-potential $\mathcal{R} : \U \to \R$ given by 
$$
   \mathcal{D} ( v) = \int_\Omega \eta | \nabla v |^2 \, dx = \partial_v \mathcal{R} ( v) [ v ]  ,
	\qquad    \mathcal{R} ( v) = \tfrac12\int_\Om \eta | \nabla v |^2 \, dx .
$$


It remains to define the interface energy. Following \cite{OrtizPandolfi_IJNME99} we introduce a function $\hat\psi : \R_+ \to \R_+$ which satisfies the following properties: 

\begin{itemize}
 \item[{\it (H1)}] $\hat\psi$ is concave, with $\hat\psi(0)=0$, $\hat\psi(w) >0$ for $w>0$, $\hat \psi(w)=\hat \psi(\xi_c)$ for $w\geq \xi_c >0$ ;
 \item[{\it (H2)}] $\hat\psi$ is of class $C^1$ in $[0,+\infty)$ and of class $C^2$ in $[0,\xi_c]$.
\end{itemize}

The above conditions imply that $\psi$ is non-decreasing. Moreover, let $\beta >0$ be defined as 
\begin{equation}\label{stima_beta}
- \beta = \min \{ \hat{\psi}'' (w) \text{ for } w\in[0,\xi_c] \} .
\end{equation}
Note that by the concavity and the $C^2$-regularity of $\hat\psi$ in $[0, \xi_c]$ we have $\beta \ge 0$; however, $\beta=0$ contradicts the $C^1$-regularity in $[0,+\infty)$, since $\hat\psi$ would be linear in $[0, \xi_c]$ and constant in $[\xi_c , +\infty)$.
We denote $\R^2_+ = \R \times \R_+$ and define 
$\psi : \R^2_+ \to \R_+$ by
\begin{equation} \label{e.defpsi}
	\psi ( w, \xi) = 
	\begin{cases} 
		\hat{\psi}( |w| ) 	& \text{if $ | w | \ge \xi $,}  \\
		{\displaystyle \hat\psi( \xi) - \hat{\psi}' ( \xi) \left( \frac{\xi^2 - w^2}{2\xi} \right)}   &  \text{if $| w | < \xi$,}
	\end{cases}
\end{equation}
(see Figure \ref{fig1} for an example). The behaviour for $| w| < \xi$ (and $\xi >0$) corresponds to an elastic unloading. 
Then, the interface cohesive energy takes the form 
$$
	 \Psi ( u , \xi) = \int_K \psi ( \jumpo{u} , \xi ) \, d\mathcal{H}^1  .
$$
We remark that for technical reasons $\psi$ is defined in the whole $\R^2_+$, even if the evolution will  take place in the cone $\{ (w,\xi) \in \mathbb{R}^2_+ : | w | \le \xi \}$. 

In order to better understand the thermodynamics of the system, it is interesting to split $\psi$ 
in terms of stored and dissipated densities, i.e.~$\psi ( w,\xi) = \psi_s ( w,\xi)   + \psi_d (\xi) $ where $\psi_d (\xi) = \psi( 0 , \xi)$. Accordingly, we introduce the stored and dissipated interface energies
\begin{equation}\label{splitting}
      \Psi_s ( u , \xi ) = \int_K \psi_s ( \jumpo{u} , \xi ) \, d\mathcal{H}^1 ,\qquad  \Psi_d ( \xi ) = \int_K \psi_d ( \xi ) \, d\mathcal{H}^1 . 
\end{equation}
We will show in Lemma \ref{l.psi} below that the functional $u\mapsto \Psi(u,\xi)$ is $\lambda$-convex; this means that there exists $\lambda<0$ such that 
\begin{align}\label{lambda_convexity}
 \Psi(u,\xi)-\lambda\|\jumpo{u} \|^2_{L^2(K)},
\end{align}
is convex for all $\xi\in\Xi$.

The above assumptions are  sufficient to prove existence and energy identity of weak solutions; in order to prove existence of strong solutions, improving the time regularity, we will assume the following additional hypotheses:
\begin{itemize}
\item[{\it (H3)}] $\hat\psi'$ is concave on  $[0,\xi_c]$;
\item[{\it (H4)}] there exists $c>0$ such that
\begin{equation} \label{e.unicoe}
	 \mu \| u \|^2_{\U} - \beta \, \| \jumpo{u} \|^2_{L^2(K)} \ge c \, \| u \|^2_{\U}  . 
\end{equation}
\end{itemize}

Note that under assumption {\it (H4)}  the functional 
\begin{align}\label{strict_convexity}
 u\mapsto \tfrac12 \mu \| u \|^2_{\U} +\Psi(u,\xi)
\end{align}
is strictly convex for all $\xi\in \Xi$. Indeed by \eqref{stima_beta} the map $w \mapsto \psi(w,\xi) + \tfrac12 \beta w^2$ is convex and thus, by linearity of the trace, the functional $u \mapsto \Psi ( u , \xi) + \frac12 \beta \| \jumpo{u} \|^2_{L^2}$ is convex. Writing
\begin{align*}
	\tfrac12 \mu \| u \|^2_{\U} +\Psi(u,\xi) = \big( \tfrac12 \mu \| u \|^2_{\U} - \tfrac12 \beta \| \jumpo{u} \|^2_{L^2} \big) + \left( \Psi(u,\xi) + \tfrac12 \beta \| \jumpo{u} \|^2_{L^2} \right)
\end{align*}
we get \eqref{strict_convexity} since the first term is quadratic and positive (for $u\neq0$) and thus strictly convex.
%

\begin{remark} \normalfont By trace and Poincar\'e inequalities we have
\begin{equation}\label{trace_ineq}
  \| u \|^2_{\U}\geq \tilde c \, \|\jumpo{u}\|_{L^2(K)}^2,
\end{equation}
and thus \eqref{e.unicoe} holds for $\mu \tilde{c} > \beta$.
Alternatively, given $\mu$ and $\beta$, \eqref{e.unicoe} holds for ``small domains''. Indeed, consider the sets $l \Omega$ for $l >0$ and the functions $u_l (x) = u ( x / l )$ we have
\begin{equation} \label{e.sca}
    \|   \nabla u_l  \|_{L^2(l \Omega)}^2 = \| \nabla u \|^2_{L^2 (\Omega)} = \| u \|^2_{\U} \ge \tilde c \, \|\jumpo{u}\|_{L^2(K)}^2 =  ( \tilde c / l )  \, \|\jumpo{u_l} \|_{L^2(l K)}^2 .
\end{equation}
If $l$ is sufficiently small \eqref{e.unicoe} holds. 
\end{remark}

\subsection{An example}\label{example}

This setting complies for the simple case where $\Omega^+ = (0,L) \times (0,1)$, $\Omega^- = (0,L) \times (-1,0)$, $K = (0,L) \times \{ 0 \}$ and $\partial_D \Omega = (0,L) \times \{-1,1\}$.

\begin{figure}[!h] \begin{center} \hspace{-7pt} \includegraphics[scale=1]{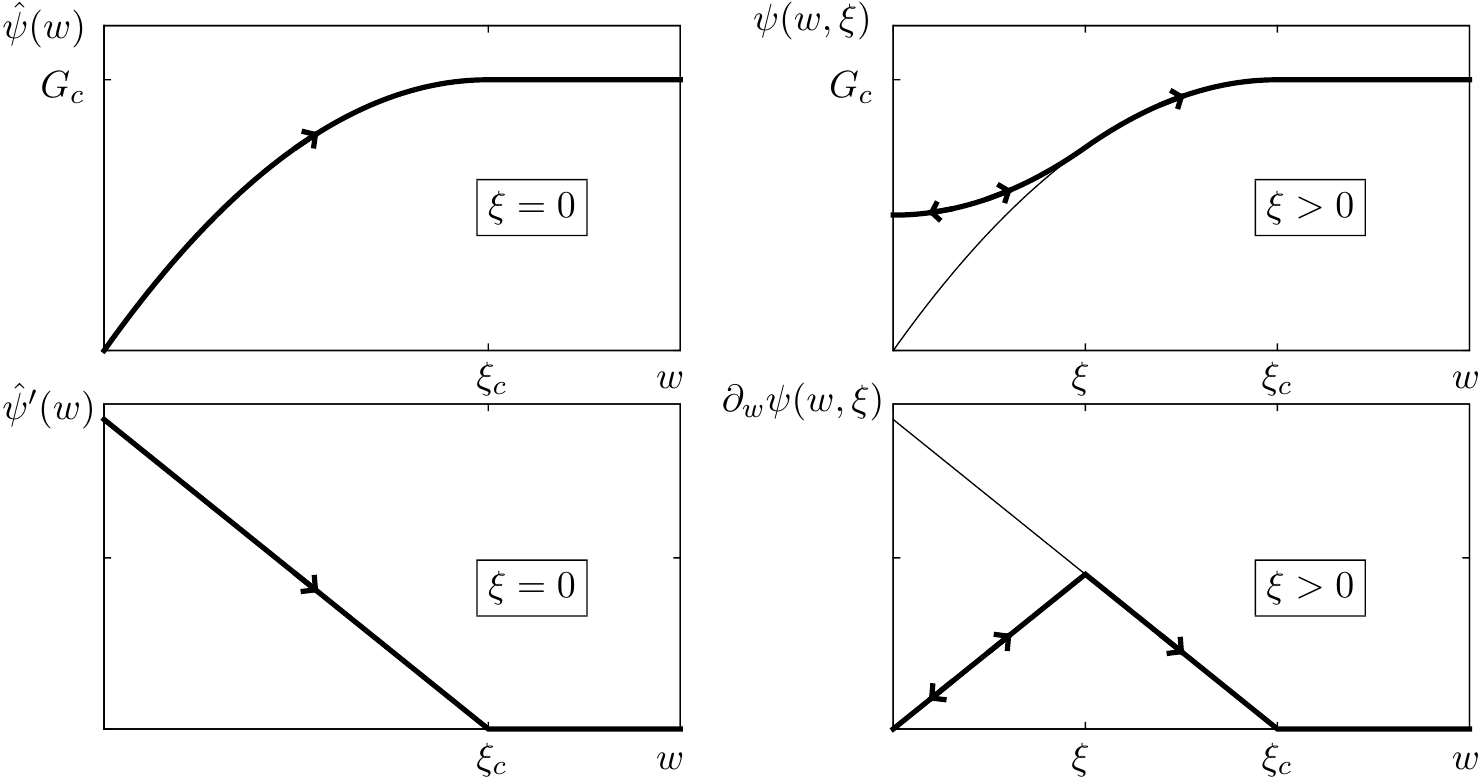} 
\end{center} \caption{\label{fig1} An example of cohesive potential (top) and the associated traction-separation law (bottom).} 
\end{figure}

We provide here a prototype example of cohesive potential (see Figure \ref{fig1}). For $G_c, \xi_c >0$ define $\hat\psi  :  [0,+\infty) \to [0,+\infty)$ as 
\begin{equation} \label{e.psi-es}
	\hat\psi ( w ) = \begin{cases}
		G_c (w/\xi_c) [ 2  - (w/\xi_c) ] & 0 \le w \le \xi_c \\[3pt]
		G_c & w > \xi_c .
	\end{cases}
\end{equation}
Note that $\hat\psi$ is of class $C^1$ and that 
$$
	\hat{\psi}'( w ) = \begin{cases}
	2 ( G_c / \xi_c) (1 - (w/\xi_c) ) & 0 \le w  \le \xi \\[3pt]
	0 & w > \xi_c.
	\end{cases}
$$
Then, for $\xi >0$, the potential $\psi$ takes the form 
\begin{align}\label{ex}
	\psi ( w , \xi ) = \begin{cases}
	\hat\psi(|w|) & |w |> \xi , \\[3pt]
	\tfrac12 ( \hat{\psi}'(\xi) / \xi  ) w^2 + \big( \hat\psi(\xi) - \tfrac12  \hat{\psi}'(\xi) \xi \big) &  | w| \le \xi .
	\end{cases} 
\end{align}
Note that
$$
 \psi_d (\xi) =  \begin{cases} 
 		\hat\psi(\xi) - \tfrac12  \hat{\psi}'(\xi) \xi = (G_c / \xi_c) \, \xi ,  &  0 \le \xi \le \xi_c \\
 		G_c & \xi \ge \xi_c .
 				\end{cases}
$$
while
$$        \psi_{s} (w,\xi) = \tfrac12 ( \hat{\psi}'(\xi) / \xi  ) \, w^2 = \tfrac12 c_\xi w^2 \quad \text{for $|w | \le \xi$ and $\xi>0$} , $$
where we have set 
$$  c_\xi:=\frac{\hat{\psi}'(\xi)}{\xi}.  $$
In particular, 
$\psi_s ( \cdot ,\xi) $ is quadratic for $| w| \le \xi$ and $\xi >0$; hence, the elastic domain depends on $\xi$ and vanishes for $\xi =0$, while $c_\xi$ becomes singular when $\xi \to 0^+$,  because $\hat\psi'(0) > 0$. Finally, note that in the energy density the singularity of $c_\xi$ is in some sense ``balanced'' by $w$, indeed, for $| w | \le \xi$ it holds $c_\xi w ^2 \le c_\xi \xi^2 \le \hat{\psi}' (\xi)  \xi \le C \xi$.

{\begin{remark} \normalfont
Assume $\hat{\psi}$ of the form \eqref{e.defpsi}. If, for some $\alpha >0$, $\psi_d$ is of the form
$$
\psi_d(\xi) = \begin{cases}
\alpha  \, \xi &  0\le \xi \le \xi_c \\
\alpha \xi_c  &   \xi > \xi_c ,
\end{cases}
$$
then, solving the ODE $ \hat\psi(\xi) - \tfrac12  \hat{\psi}'(\xi) \xi = \alpha  \, \xi$ (which defines $\psi_d$)  it turns out that $\hat\psi$ is of the form \eqref{e.psi-es}. 
\end{remark}}

\begin{figure}[!h] \begin{center} \includegraphics[scale=1]{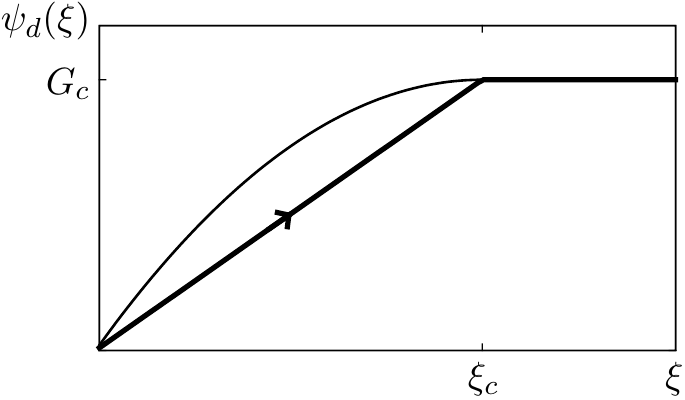} \caption{\scriptsize \label{f.du} Dissipated energy $\psi_d$ (bold) compared with the energy $\hat\psi$.} \end{center} \end{figure}

\subsection{Directional derivatives and subdifferential} 

We will denote by $\partial_w \psi ( w, \xi ; \phi)$ the directional derivative of $\psi(\cdot,\xi)$ with respect to $\phi$, and by $\partial_w \psi ( w, \xi )$ the subdifferential. We will also denote by  $\partial_w \psi ( w, \xi )$ the partial derivative, where it exists, and we will further note by $\partial^\pm_w \psi ( w, \xi )$ the left and right derivatives.
To compute explicitly the derivatives of $\psi$, it is convenient to consider separately the interior, the exterior, and the boundary of the cone $\{ | w | \le \xi \}$ in $\R^2_+$. In the set $\{ |w| < \xi \}$ (for $\xi >0$) we have 
\begin{equation} \label{e.stella}
	\partial_w \psi (w, \xi) = \bigg( \frac{\hat{\psi}' ( \xi)}{\xi} \bigg) w .
\end{equation}
In the set $\{ | w| > \xi \}$ we have $\psi ( w, \xi) = \hat\psi ( | w | )$, hence 
$$
	\partial_w \psi (w, \xi) =  \hat{\psi}' ( | w| ) \, \mathrm{sign}(w)  .
$$
It is easy to check that the above derivatives coincide in the set $\{ | w | = \xi \neq 0 \}$; more precisely, for $w = \xi >0$ we have
$$
    \partial^+_w \psi ( w, \xi ) = \hat{\psi}' ( \xi )  = \partial^-_w \psi ( w, \xi) ,
$$
and similarly for $w = -\xi <0$. As a consequence $\partial_w \psi ( w, \xi)$ is continuous in $\R^2_+ \setminus (0,0)$.
On the contrary, the density $\psi$ is not differentiable in the origin where we have only the directional derivatives 
\begin{align}\label{der_0}
	\partial_w \psi ( 0 , 0 ; \phi) = \lim_{h \to 0^+}  \frac{\psi ( h \phi , 0) }{h} = \hat{\psi}'(0) | \phi | 
	\quad \text{for every $\phi$.} 
\end{align}
It is important to observe that directional derivatives are positively $1$-homogeneous, i.e., $\partial_w \psi ( w , \xi ; \kappa \phi) = \kappa \partial_w \psi ( w , \xi ; \phi)$ for every $\phi \in \R$ and $\kappa >0$ and that in general we have
\begin{equation}\label{e.derctrl} | \partial_w \psi ( w , \xi ; \phi ) | \le \hat\psi' ( 0) \, | \phi | . 
\end{equation}

In a similar way we can compute the directional derivatives $\partial_\xi \psi ( w, \xi ; \zeta)$ and the partial derivative $\partial_\xi \psi ( w, \xi)$, where it exists. In the set $\{ |w| < \xi \}$ (for $\xi >0$) we further distinguish the cases: $\xi<\xi_c$, $\xi= \xi_c$, and $\xi> \xi_c$. 
In the former case we have
\begin{equation}\label{derivative_xi}
	\partial_\xi \psi (w,\xi) = - \tfrac12 \big( \hat\psi'' ( \xi) \xi - \hat\psi' (\xi) \big) \bigg( \frac{\xi^2 - w^2}{\xi^2} \bigg)  = - \bigg( \frac{\hat{\psi}' (\xi)}{2\xi} \bigg)' (\xi^2 - w^2)  . 
\end{equation}
For $\xi > \xi_c$ we have $\psi ( \xi , w) = \hat\psi( \xi_c)$ and thus $\partial_\xi \psi ( \xi , w) = 0$.
For $\xi= \xi_c$ we have only directional derivatives, more precisely
$$
   \partial^-_\xi \psi (w,\xi_c) = - \tfrac12 \big( \hat\psi'' ( \xi_c) \xi_c - \hat\psi' (\xi_c) \big) \bigg( \frac{\xi_c^2 - w^2}{\xi_c^2} \bigg) , \qquad
   	 \partial^+_\xi \psi (w,\xi_c) = 0 .
$$ 
In the set $\{ | w| > \xi \}$ we simply have $\partial_\xi \psi (w,\xi) =  0$ because $\psi ( w, \xi) = \hat\psi ( | w | )$. It is not difficult to check that in the set $\{|w| = \xi>0\}$ we have 
$$
    \partial_\xi^+ \psi ( w ,\xi ) = \partial^-_\xi \psi ( w, \xi) = 0.
$$
Finally, in the origin, by \eqref{ex},  we have
$$
    \partial_\xi \psi ( 0 , 0 ; \zeta) = \lim_{h \to 0^+} \frac{\psi ( 0 , h \zeta ) }{h} = 
	\lim_{h \to 0^+} \frac{1}{h} \big( \hat{\psi}(h \zeta ) - \hat{\psi}(0)  - \tfrac12 \hat\psi'(h \zeta) h \zeta \big) = \tfrac12 \hat\psi'(0) \zeta \quad \text{for every $\zeta \ge 0$.}
$$ 
 {Note that there exists $C>0$ such that 
 $$		|\partial_\xi \psi ( w , \xi ; \zeta )  | \le C | \zeta |  $$ 
 for every $(w,\xi) \in \R^2_+$ and every admissible $\zeta$.}

\begin{lemma} \label{l.psi} Under assumption (H1) and (H2) the density $\psi$ satisfies the following properties: 
\begin{itemize}
\item[$\bullet$] $\psi$ is Lipschitz continuous in $\mathbb{R}^2_+$ {and differentiable in $\{ 0 < | w | < \xi_c \}$};
\item[$\bullet$] $\partial_w \psi$ is continuous in $\mathbb{R}^2_+ \setminus (0,0)$;
\item[$\bullet$] $\psi ( w, \cdot)$ is monotone non-decreasing;
\item[$\bullet$] $\psi ( \cdot, \xi )$ is pair, monotone non-decreasing in $\R_+$ and $\lambda$-convex in $\R$, uniformly  w.r.t.~$\xi \in \R_+$.
\end{itemize}
\end{lemma}

\proof By the above computations it is easy to see that: the directional derivatives are bounded in the whole $\R^2_+$, that $\partial_w \psi$ is continuous in $\R^2_+ \setminus (0,0)$, and that $\partial_\xi \psi$ is continuous in $\{ 0 < | w | < \xi _c \}$. 
Moreover, by concavity and monotonicity, $ \hat\psi'(0) \ge \hat\psi'(\xi_1) \ge \hat\psi' ( \xi_2) \ge 0$ for every $0 \le \xi_1 \le \xi_2$; it follows that $\xi \mapsto \hat{\psi}' ( \xi) /\xi $ is monotone non-increasing, for $\xi >0$; hence $\partial_\xi \psi ( w, \xi) \ge 0$ and thus $\psi ( w , \cdot) $ is monotone non-decreasing. It is easy to check that  $\psi ( \cdot, \xi )$ is pair and monotone non-decreasing in $\R_+$. By \eqref{stima_beta} and \cite[Lemma 2.4]{NScala_NARWA17} we get the uniform $\lambda$-convexity. \qed

\begin{remark} By Lemma \ref{l.psi} the dissipated energy $\Psi_d$ is monotone non-decreasing. Moreover the energy $\Psi ( \cdot, \xi)$ is $\lambda$-convex in $L^2(K)$, uniformly with respect to $\xi$. \end{remark}

\begin{lemma} \label{l.usc-psi'} 
Assume (H1) and (H2).
If $(w_n, \xi_n) \to (w,\xi)$ in $\R^2_+$ then $\limsup_{n \to +\infty} \partial_w \psi ( w_n , \xi_n ; \phi) \le \partial_w \psi ( w, \xi ; \phi)$. If $\phi_n \to \phi$ then $\lim_{n \to +\infty} \partial_w \psi ( w_ , \xi ; \phi_n) = \partial_w \psi ( w, \xi ; \phi)$. 

\end{lemma}

\proof 
If $(w,\xi) \neq (0,0)$ the statement is true because $\partial_w \psi ( w , \xi ; \phi) = \partial_w \psi ( w , \xi ) \phi$ with $\partial_w \psi$ continuous away from the origin. If $(w,\xi) = (0,0)$ the same conclusion holds since  $\partial_w \psi ( w_n ,\xi_n ; \phi) \le \hat\psi' (0) | \phi |  =  \partial_w \psi ( 0, 0 ; \phi)$. The convergence with respect to $\phi_n$ is obvious, considering again the cases $(w,\xi) = (0,0)$ and $(w,\xi) \neq (0,0)$. \qed 


\medskip
Under the above assumptions we can introduce $\partial_u \F ( t, u , \xi ; \phi )$,  the directional derivative of the energy $\F ( t, \cdot, \xi)$ in the direction $\phi$, which reads 
$$\partial_u \F ( t, u , \xi ; \phi ) = d \E (  u ) [ \phi ] - ( f(t) , \phi )_\U + \partial_u \Psi ( u , \xi ; \phi ), $$
where
$$
	d \E (u) [\phi] = \int_\Omega \mu \nabla u \cdot \nabla \phi \, dx , 
\quad \partial_u \Psi ( u , \xi ; \phi )  = \int_K \partial_w \psi ( \jumpo{u} , \xi ; \jumpo{\phi}  ) \, d\mathcal{H}^1 .
$$
%
%
Moreover, we recall that $\zeta \in \U^*$ is a (Clarke directional) subderivative for $\F$ in $(t,u,\xi)$ if 
$$
	\partial_u \F ( t, u , \xi ; \phi ) \ge ( \zeta , \phi)_\U \quad \text{for every $\phi \in \U$},
$$
and that $\zeta \in \U^*$ belongs to the (limiting) subdifferential $\partial_u \F  ( t, u , \xi )$ if
$$
	\liminf_{w \to u} \frac{\F ( t, w , \xi) - \F (t,u,\xi) - ( \zeta , w- u)_\U }{\| w - u\|_{\U} } \ge 0 \,.
$$	
Here the convergence $w\rightarrow u$ is intended with respect to the strong topology of $\U$. 
Thanks to the properties of $\F$ these two notions coincide, we have indeed the following statement.

\begin{lemma} \label{l.sub=var} Assume (H1) and (H2) and let $(t,u,\xi) \in [0,T] \times \U \times \xi$ and $\zeta \in \U^*$; $\zeta$ is a subderivative for $\F$ in $(t,u,\xi)$ if and only if it belongs to the subdifferential $\partial_u \F (t,u,\xi)$. In particular the equilibrium condition $\partial_u \F  ( t, u , \xi ; \phi ) \ge 0$  for every $\phi \in \U$ is equivalent to the inclusion $\partial_u \F  ( t, u , \xi ) \ni 0$.
\end{lemma}

\proof If $\zeta \in \partial_u \F  ( t, u , \xi )$  and $\phi \in \U \setminus \{ 0 \}$, the change of variable $w = u + h \phi$, for $h >0$ leads to 
$$
	\liminf_{h \to 0^+} \frac{\F ( t, u + h \phi ,\xi) - \F ( t, u, \xi)}{h \| \phi \|_\U} \ge \frac{ ( \zeta , \phi)_\U}{ \| \phi \|_\U} ,
$$
which gives $\partial_u \F  ( t, u , \xi ; \phi ) \ge ( \zeta , \phi )_\U$.

On the other hand, by $\lambda$-convexity the auxiliary functional 
$$	\F_{\lambda, u} ( t, w , \xi) = \F ( t, w , \xi) - \lambda \| w - u \|^2_\U   $$
is convex  in $w$. Note that $\F_{\lambda, u} ( t, u , \xi) = \F (t,u,\xi)$, 
$\partial_u \F_{\lambda, u} ( t, u , \xi; \phi) = \partial_u \F (t,u,\xi; \phi)$ for every $\phi \in \U$ and that 
$$
	\liminf_{w \to u} \frac{\F_{\lambda,u} ( t, w , \xi) - \F_{\lambda,u} (t,u,\xi) - ( \zeta , w- u)_\U }{\| w - u\|_{\U} } = 
	\liminf_{w \to u} \frac{\F ( t, w , \xi) - \F (t,u,\xi) - ( \zeta , w- u)_\U }{\| w - u\|_{\U} } .
$$
By convexity we have 
$$\F_{\lambda,u} ( t, w , \xi) - \F_{\lambda,u} ( t, u , \xi) \ge \partial_u \F_{\lambda,u} (t,u,\xi; w-u) = \partial_u \F  (t,u,\xi; w-u).$$ 
Hence, 
$$
	\liminf_{w \to u} \frac{\F_{\lambda,u} ( t, w , \xi) - \F_{\lambda,u} (t,u,\xi) - ( \zeta , w- u)_\U }{\| w - u\|_{\U} } \ge 
	\liminf_{w \to u} \frac{\partial_u \F (t,u,\xi; w-u) - ( \zeta , w- u)_\U }{\| w - u\|_{\U} } .
$$
Let us consider a sequence $w_k$ such that 
\begin{align*}
\liminf_{w \to u} \frac{\partial_u \F (t,u,\xi; w-u) - ( \zeta , w- u)_\U }{\| w - u\|_{\U} } & = 
	\lim_{k \to \infty} \frac{\partial_u \F (t,u,\xi; w_k-u) - ( \zeta , w_k - u)_\U }{\| w_k - u\|_{\U} } \\
	& = \lim_{k \to \infty} \partial_u \F (t,u,\xi; \phi_k ) - ( \zeta , \phi_k)_\U ,
\end{align*}
where $\phi_k = (w_k - u) / \| w_k - u \|_\U $. We extract a further subsequence, non relabelled, such that $\phi_k$ converge weakly to $\phi$ in $\U$ and (by compactness of the trace) strongly in $L^2 (K)$. Then, being $\partial_u \F ( t, u , \xi ; \phi ) = d \E (  u ) [ \phi ] - ( f(t) , \phi )_\U + \partial_u \Psi ( u , \xi ; \phi )$, by Lemma \ref{l.psi} and Lemma \ref{l.usc-psi'} we have $\partial_u \Psi ( u , \xi ; \phi_k ) \to \partial_u \Psi ( u , \xi ; \phi )$ and then if $\zeta$ is a subderivative we have 
$$
	\lim_{k \to \infty} \partial_u \F (t,u,\xi; \phi_k ) - ( \zeta , \phi_k)_\U  = \partial_u \F (t,u,\xi; \phi) -  ( \zeta , \phi )_\U  \ge 0 ,
$$
which concludes the proof. \qed

\section{Main results}\label{mainresults_sec}

\subsection{Weak solutions and energy identity\label{genset}}

We introduce the concept of weak solution we are going to investigate in what follows:

\begin{definition} \label{d.wsol} Given $(u_0 , \xi_0) \in \U \times \Xi$, with $| \jumpo{u_0} | \le \xi_0$, and  $v_0 \in L^2(\Om)$, a couple $(u,\xi) \in W^{1,2} (0,T; \U \times \Xi)$ with $\ddot u \in L^2 (0,T ; \U^*)$ is a weak solution 
with initial conditions $u_0$ and $v_0$ if 
\begin{align}\label{maineq}
	\begin{cases}
	\rho \ddot u (t) + \partial_u \F ( t , u (t) , \xi(t) ) + \partial_v \mathcal{R} ( \dot{u})  \ni0 ,  & \text{for a.e.~$t \in (0,T)$}, \\[5pt]
	\dot{\xi} (t)  ( \xi(t) -  | \jumpo{u(t)} | ) = 0 \ \text{ and } \ | \jumpo{u(t)} | \le \xi (t) ,   & \text{for a.e.~$t \in (0,T)$,} \\[4pt]
	u (0) = u_0 , \ \dot{u} (0) = v_0 .
	\end{cases}
\end{align}
By Lemma \ref{l.sub=var} the above differential inclusion  in variational form reads 
 $$
( \rho \ddot u (t) , \phi )_\U + \partial_u \F ( t , u (t) , \xi(t)  ; \phi ) + \langle \eta\nabla \dot u (t) , \nabla \phi \rangle \ge 0 ,  \ \text{for every $\phi \in \U$.} 
$$
%
\end{definition}

\noindent In the next sections we will prove the following result. 

\begin{theorem} \label{t.teor} Assume that hypotheses (H1) and (H2) hold and that $f\in W^{1,2} (0,T;\U^*)$. Let $(u_0 , \xi_0) \in \U \times \Xi$, with $| \jumpo{u_0} | \le \xi_0$, and $v_0 \in L^2(\Om)$. Then, there exists a weak solution $(u,\xi)$, in the sense of Definition \ref{d.wsol}, 
which satisfies the energy identity 
\begin{align}
	\E ( u( t^*) ) + \Psi( u(t^*) , \xi (t^*)) + \K ( \dot{u} (t^*) )  & =  \, \E ( u_0 ) + \Psi ( u_0, \xi_0) + \K ( v_0 ) \nonumber \\ & \quad + \int_0^{t^*} \mathcal{P}_{ext} ( t , \dot{u} (t) ) \, dt - \int_0^{t^*} \mathcal{D} ( \dot{u} (t) )  \,dt ,  \label{e.bilen}
\end{align}
for every $t^* \in [0,T]$.
\end{theorem}

Splitting the stored and dissipated parts of $\Psi$, as in \eqref{splitting}, the above energy identity can be written also in the following way:
\begin{align}
	\E ( u( t^*) ) + \Psi_s ( u(t^*) , \xi (t^*)) + \K ( \dot{u} (t^*) )  & =  \, \E ( u_0 ) + \Psi_s ( u_0, \xi_0) + \K ( v_0 ) + \int_0^{t^*} \mathcal{P}_{ext} ( t , \dot{u} (t) ) \, dt  \nonumber \\ & \quad - \int_0^{t^*} \mathcal{D} ( \dot{u} (t) ) \, dt - \int_0^{t^*} \partial_\xi \Psi_d ( \xi(t) ; \dot{\xi}(t)) \, dt  , \label{e.bilen.bis}
\end{align}
where $\partial_\xi \Psi_d ( \xi(t) ; \dot{\xi}(t))$ is the dissipation on the interface.

\subsection{Solutions with higher time regularity}

Under stronger assumptions on the data, we are able to refine Theorem \ref{t.teor} as follows.

\begin{theorem}\label{d.wsol2} 
 Besides the hypotheses of Theorem \ref{t.teor}, assume that (H3) and (H4) hold and that $f\in W^{1,2}(0,T;L^2(\Om))$. Moreover, assume that $v_0 \in \U$ with $\jumpo{v_0}=0$ and that there exists 
$w_0 \in L^2(\Omega)$ such that 
\begin{equation}	\label{e.w0}
\rho w_0  + \partial_u \F ( t , u_0 (t) , \xi_0 (t) ) + \partial_v \mathcal{R} ( {v}_0)  \ni0 .  
\end{equation}
Then, there exists a weak solution $(u,\xi)$, in the sense of Definition \ref{d.wsol}, satisfying the energy balance \eqref{e.bilen} and the further regularity
 \begin{align}\label{regularity}
  &u\in W^{2,\infty}(0,T;L^2(\Om))\cap W^{1,\infty}(0,T;\U).
 \end{align}
\end{theorem}

\begin{remark}\label{equivalence_initialdata} \normalfont
The function $w_0$ appearing in Theorem \ref{d.wsol2} plays the role of $\ddot{u}_0$, thus \eqref{e.w0} means that the initial acceleration belongs to $L^2$. 
 Again by Lemma \ref{l.sub=var},  condition \eqref{e.w0} is equivalent to 
 \begin{align}\label{e.w0equiv}
  \langle \rho w_0,\phi\rangle +\partial_u \F ( t , u_0 , \xi_0  ; \phi ) + \langle \eta\nabla v_0 (t) , \nabla \phi \rangle \ge 0 ,  \ \text{for every $\phi \in \U$.}
 \end{align}
This is readily satisfied if   
 \begin{align}\label{u_0minimo}
  u_0\in \argmin\{\mathcal{F}  (0, u,\xi_0) +\langle\eta\nabla v_0,\nabla u\rangle+\langle \rho w_0,u\rangle, \;u\in \U \}.
 \end{align}
Notice that such $u_0$ always exists and is also unique thanks to \eqref{strict_convexity}. In particular \eqref{e.w0} and \eqref{u_0minimo} turn out to be equivalent.
\end{remark}

\subsection{Strong solutions \label{s.poly}}

In this section we consider $\Omega^\pm$ to be (curvilinear) polygonal domains, in the sense of \cite{Grisvard}. In this setting, under the assumptions of Theorem \ref{d.wsol2} the weak solutions given by Theorem \ref{d.wsol2} satisfy also the system of PDEs \eqref{e.pde}.

\begin{theorem}\label{teo_curvilinear} Besides the hypotheses of Theorem \ref{d.wsol2} assume that $\Om^\pm$ are (curvilinear) polygonal domains. Let $(u,\xi)$ be a solution provided by Theorem \ref{d.wsol2} and denote  by $\sigma (t) = \mu \nabla u (t) + \eta \nabla \dot{u} (t)$ the visco-elastic stress. Then, the differential inclusion in $\U^*$
$$
	\rho \ddot u (t) + \partial_u \F ( t , u (t) , \xi(t) ) + \partial_v \mathcal{R} ( \dot{u})  \ni0 , 
$$
is equivalent to the following system of partial differential equations: 
\begin{gather}  \label{e.pde}
\begin{cases}
 \rho \ddot{u}(t) - \Div \, \sigma(t)  = f (t) \qquad& \hbox{in $\Omega$,}\cr
 u(t)=\dot u(t) =0 &\hbox{in $\partial_D\Omega^\pm$,}\cr
\sigma (t) \nu =0 &\hbox{in $\partial_N \Omega^\pm$, }\cr
 \sigma^+(t)\nu = \sigma^-(t)\nu  \in \partial_w \psi \big( \jumpo{u (t)},\xi(t) \big) &\hbox{in $K$.}  
\end{cases}
\end{gather}

In particular $\sigma^\pm(t) \nu \in L^\infty (K)$ and $\sigma^+ (t) \nu = \partial_w \psi ( \jumpo{u (t)} ,  \xi(t) )$ if $\xi(t)>0$ while $| \sigma^+(t) \nu | \le \hat{\psi}' (0)$  if $\xi(t)=0$. 

\end{theorem}
\begin{remark} \normalfont
Let us spend some words on the hypothesis \eqref{e.w0}. 
 In this more regular setting where $\Om^\pm$ are (curvilinear) polygonal domains, we might write condition \eqref{e.w0} in a more treatable way.
 
 Precisely, thanks to the regularity of the domain, we can now integrate by parts equation \eqref{e.w0} (or \eqref{e.w0equiv}) and obtain that this condition is equivalent to the system
 \begin{gather}  \label{e.sigma0}
\begin{cases}
 \rho w_0 - \Div \, \sigma_0  = f(0) \qquad& \hbox{in $\Omega$,}\cr
 u_0=v_0 =0 &\hbox{in $\partial_D\Omega^\pm$,}\cr
\sigma_0 \nu =0 &\hbox{in $\partial_N \Omega^\pm$, }\cr
 \sigma^+_0\nu = \sigma^-_0\nu  \in \partial_w \psi \big( \jumpo{u_0},\xi_0 \big) &\hbox{in $K$,}  
\end{cases}
\end{gather}
 where we have set $\sigma_0=\mu\nabla u_0+\eta\nabla v_0$. The first condition also reads
 \begin{align}\label{reg_w0}
  \rho w_0=\mu\Delta u_0+\eta\Delta v_0+f(0).
 \end{align}
Since $f(0)\in L^2(\Omega)$ the existence of $w_0$ satisfying \eqref{reg_w0} is achieved as soon as 
$$\mu\Delta u_0+\eta\Delta v_0\in L^2(\Omega).$$
This is a natural requirement, since this condition is always satisfied by $(u,v)$ during the evolution.
The first equation in \eqref{e.sigma0} can then be seen as a compatibility condition for the initial stress $\sigma_0$, since it just requires that $-\Div \sigma_0\in L^2(\Omega)$. Also the other conditions in \eqref{e.sigma0} are natural and can be viewed as compatibility conditions for $\sigma_0$.

The additional hypothesis  $\jumpo{v_0}=0$ has the following interpretation. At the points  where $\xi_0=0$ we necessarily have $ \jumpo{\dot u(0)}=0$, indeed, by definition of $\xi$, we would have $\xi(t)=0$ and $|\jumpo{u(t)}|=0$ for $t<0$. Therefore $\jumpo{v_0}=0$ is a natural condition on the set $\{\xi_0=0\}\subset K$. To simplify the arguments we assume that $ \jumpo{v_0}=0$ on the whole $K$. This condition is possibly generalizable but  we prefer not to weaken it for technical reasons, that will be evident  in the proof of Proposition \ref{p.kkomp} below, and in order to provide a more clear exposition.
 
\end{remark}
In view of the previous remark, we can state the equivalent of Theorem \ref{d.wsol2}.

\begin{theorem}
Let us assume the hypotheses of Theorem \ref{d.wsol}, and suppose  $\Om^\pm$ are (curvilinear) polygonal domains. Assume also that (H3) and (H4) hold and $f\in W^{1,2}(0,T;L^2(\Om))$. If $u_0\in \U$, $v_0\in \U$, $\xi_0\in \Xi$ satisfy
\begin{gather}  \label{e.initialdata}
\begin{cases}
\sigma_0 \nu =0 &\hbox{in $\partial_N \Omega^\pm$, }\cr
 \sigma^+_0\nu = \sigma^-_0\nu  \in \partial_w \psi \big( \jumpo{u_0},\xi_0 \big) &\hbox{in $K$,} \cr
 \jumpo{v_0}=0&\hbox{in $K$,}\cr
 {|\jumpo{u_0}|} \leq\xi_0&\hbox{in $K$,}
\end{cases}
\end{gather}
and $$\Div\sigma_0\in L^2(\Omega),$$
then the same conclusions of Theorems \ref{d.wsol2} and  \ref{teo_curvilinear} hold.
\end{theorem}

\section{A preliminary regularized evolution}\label{preliminaryevolution}

For sake of simplicity we will assume that $\rho$, $\mu$ and $\eta$ are constant, i.e., $\rho^+ = \rho^-$ etc. We also fix a positive number $\bar \xi$.
%
%


%
%

\subsection{Time discrete evolution}

Assume that the cohesive potential $\Psi$ enjoys conditions (H1) {and (H2)}. Moreover, assume that  $\xi_0 \ge \bar\xi > 0$ a.e.~on $K$; note that under the latter assumption the energy functional turns out to be differentiable since the cohesive potential $\psi ( \jumpo{u} , \xi)$ is singular only when $ \jumpo{u}= \xi =0$. Thus, for $\xi \ge \bar\xi$, the functional $u \mapsto  \Psi ( u ,  \xi) $ is Fr\'echet differentiable and 
\begin{equation}  \label{e.Psiceps}
	\partial_u \Psi  ( u , \xi ) [ \phi ]  = \int_K \partial_w \psi ( \jumpo{u} , \xi) \jumpo{\phi} \, d\mathcal{H}^1 = \int_K c_{\xi}\jumpo{u} \jumpo{\phi} \, d\mathcal{H}^1 , 	\qquad 	c_\xi =  \hat{\psi}'(\xi) / \xi ,
\end{equation}
if $| \jumpo{u} | \le \xi$, see \eqref{e.stella}. 

Let $\tau_n = T/n$ and let $t_{n,k} = k \tau_n$ for $k =-1, \ldots, n$. Since we assume $f\in W^{1,2}(0,T;\U^*)$, for $k=0, \ldots, n$ we set, by continuity,
\begin{equation*}
f_{n,k}=f(t_{n,k}).
\end{equation*}
We denote by $f_{n}^\sharp$ the piecewise-constant, left-continuous interpolant of $f_{n,k}$, i.e.
\begin{align*}
 &f^\sharp_n(t):=f_{n,k} \quad \text{ for }t\in(t_{n,k-1},t_{n,k}].
\end{align*}
We also introduce the piecewise-affine interpolant of $f_{n,k}$, namely
\begin{align}
 &f_n(t):=f_{n,k-1}+(t-t_{n,k-1}) \tau_n^{-1} (f_{n,k}-f_{n,k-1}), \quad \text{ for }t\in[t_{n,k-1},t_{n,k}).\label{fenne}
\end{align}
Accordingly define
$$
\F_{n} (t , u , \xi) = \tfrac12  \int_\Omega \mu | \nabla u|^2 \, dx - ( f^\sharp_{n} ( t ), u )_\U + \int_K \psi ( \jumpo{u} , \xi) \, d\mathcal{H}^{1} .
$$
We will define, by induction, two finite sequences $u_{n,k}$ and $\xi_{n,k}$, for $k=-1,...,n$. For $k=-1$, we set $u_{n,-1}  = u_0 - \tau_n v_0$.
For $k=0$, we set $u_{n,0} = u_0$ and $\xi_{n,0}= \xi_0$. 
For later convenience, we also set $\xi_{n,-1} = \xi_{n,0}$. 
For $k \ge 1$, given $u_{n,k-1}$, $u_{n,k-2}$ and $\xi_{n,k-1}$ we define 
\begin{equation} \label{e.scheme}
	\begin{cases}
		u_{n,k} \in \argmin \{ \mathcal{J}_{n,k} (  u ) : u \in \U \} , \\[3pt]
		\xi_{n,k} = \max \{   \xi_{n,k-1} , | \jumpo{u_{n,k}} | \} ,
	\end{cases}
\end{equation}
where
\begin{align} \label{e.Jnk}
	\mathcal{J}_{n,k} ( u ) & = \tfrac12 \tau_n^{-2}  \rho \|   u - 2 u_{n,k-1} + u_{n,k-2} \|_{L^2}^2 + \tfrac12 \tau_n^{-1} \eta \| \nabla (u - u_{n,k-1}) \|^2_{L^2} + \F_{n} ( t_{n,k} , u , \xi_{n,k-1} ) .
\end{align}

\begin{proposition} \label{p.J!} For $\tau_n \ll 1$ there exists a unique minimizer $u_{n,k}$ of $\mathcal{J}_{n,k} $. \end{proposition}

\proof It is enough to note that $\mathcal{J}_{n,k}$ is coercive and strictly convex in $\U$ because $\F ( t , \cdot , \xi)$ is coercive. 
 Let us check that $\mathcal{J}_{n,k}$ is strictly convex. Remember that 
 $u \mapsto \Psi ( u , \xi) - \lambda \| \jumpo{u} \|_{L^2}^2$ is convex for some $\lambda <0$. 
 Thus it is enough to check that for $\tau_n \ll 1$ the functional $u \mapsto \tau^{-1}_n \| \nabla (u  - u_{n,k-1}) \|_{L^2}^2 + \lambda \, \| \jumpo{u} \|^2_{L^2}$ is convex (recall $\lambda<0$). To this aim we write
$$
	\int_K \jumpo{u}^2 \, d\mathcal{H}^{1} = \int_K \jumpo{u - u_{n,k-1}}^2 - \jumpo{u_{n,k-1}}^2 + 2\, \jumpo{u_{n,k-1}} \, \jumpo{ u }  \, d\mathcal{H}^{1} .
$$
It is enough to see that the functional $z \mapsto \tau^{-1}_n \| \nabla z \|^2_{L^2} + \lambda \, \| \jumpo{z} \|^2_{L^2}$ is strictly convex in $\U$; indeed  it
is quadratic and non-negative (by continuity of the trace) for $\tau \ll 1$. 
\qed

\begin{remark} \label{r.sharp} \normalfont Following \cite[Lemma 5.1]{NegriVitali_IFB18} it is possible to replace (a posteriori) $\xi_{n,k-1}$ with $\xi_{n,k}$ in the incremental functional $\mathcal{J}_{n,k}$, i.e.~$u_{n,k} \in \argmin \{ \mathcal{J}^\sharp_{n,k} (  u ) : u \in \U \}$ where 
$$
	\mathcal{J}^\sharp_{n,k} ( u ) = \tfrac12 \tau_n^{-2} \rho \|   u - 2 u_{n,k-1} + u_{n,k-2} \|_{L^2}^2 + \tfrac12 \tau_n^{-1} \eta \| \nabla (u - u_{n,k-1}) \|^2_{L^2} + \F_{n} ( t_{n,k} , u , \xi_{n,k} ) 
$$
is obtained from \eqref{e.Jnk} replacing $\xi_{n,k-1}$ with $\xi_{n,k}$.
This is just a consequence of the definition of the cohesive potential, indeed $\Psi  ( u_{n,k} , \xi_{n,k-1} ) = \Psi ( u_{n,k} , \xi )  = \Psi  ( u_{n,k} , \xi_{n,k} ) $ for every $\xi_{n,k-1}\le \xi \le \xi_{n,k}$ and $\Psi ( u_{n,k} , \xi )  \ge \Psi  ( u_{n,k} , \xi_{n,k} ) $ for every $\xi \ge \xi_{n,k}$. Hence, in analogy with delamination and adhesive contact models, we can recast the update $\xi_{n,k} = \max \{ \xi_{n,k-1} , | \jumpo{u_{n,k}}  | \}$ as a minimization problem, i.e.
$$
     \xi_{n,k} \in \argmin \{  \Psi ( u_{n,k} , \xi ) : \xi \ge \xi_{n,k-1} \} = \argmin \{  \F_n (t_{n,k} , u_{n,k} , \xi ) : \xi \ge \xi_{n,k-1} \}.
$$
However, in our setting this minimization problem is not well posed, since there are infinitely many minimizers, given by any $\xi_{n,k-1} \le \xi \le \xi_{n,k}$. Among these choices, the one which naturally ensures the Karush-Kuhn-Tucker condition (in the discrete setting, see \eqref{e.ED3}, and in the limit) is $\xi_{n,k} = \max \{ \xi_{n,k-1} , | \jumpo{u_{n,k}}  | \}$, which indeed appears in \eqref{e.scheme}.
\end{remark} 

\begin{remark} \normalfont 
 \label{comparison}
The cohesive model we consider is related to adhesive models, see e.g.~\cite{RossiRoubicek_NA11,RossiThomas2017,Roubicek_SJMA13,CLO,RMP,RSZ} where, roughly speaking, the interface energy is of the form
$$
     \int_K v \jumpo{u}^2 \, d\mathcal{H}^1 , 
$$
for $0 \le v \le 1$, while the dissipation on the interface is usually of the form
$$
		\alpha \int_K | \dot{v} | \, d \mathcal{H}^1, 
$$
for some $\alpha >0$.  Here $v$ plays, in some sense, the role of our internal variable $\xi$, however, note that $v$ is monotone non-increasing in time. 
In this setting the incremental problem for $v$ reads
$$
     v_{k} \in \argmin \left\{  \int_K v \, \jumpo{u_k}^2 \, d \mathcal{H}^1 + \alpha \int_K ( v_{k-1}  - v ) \, d\mathcal{H}^1 :  v \le v_{k-1}  \right\}.
$$
By the monotonicity constraint and by convexity with respect to $v$, we have $v_k < v_{k-1}$ if and only if $ \jumpo{u_k}^2 > \alpha$. Hence the interface behaves elastically below the threshold $\alpha^{1/2}$. Therefore the mechanical behaviour is quite different from ours. 
To make it more clear, note that our energy can be written in the form 
$$
     \tfrac12\int_K c_\xi \jumpo{u}^2 \, d\mathcal{H}^1 + \int_K \psi_d ( \xi)  \, d\mathcal{H}^1 ,
$$
where $c_\xi \to +\infty$ for $\xi \to 0^+$ and in general $\psi_d$ is non-linear. From the mathematical point of view the singularity of $c_\xi$ and the non-linearity of $\psi_d$ are major sources of technical difficulties.
\end{remark}


. 

\subsection{Equilibrium and compactness}

Define $v_{n,k} = ( u_{n,k} - u_{n,k-1} ) / \tau_n$ for $k=1,..,n$, and also $\dot{f}_{n,k} = (f_{n,k} - f_{n,k-1} ) / \tau_n$.
For $k \ge 0$ we write the first term of  $\mathcal{J}_{n,k} ( u_{n,k})$ as 
$$
	\tfrac12 \tau_n^{-2}  \|   u_{n,k} - 2 u_{n,k-1} + u_{n,k-2} \|^2_{L^2} = \tfrac12 \left\|  \frac{u_{n,k} - u_{n,k-1}}{\tau_n} - \frac{ u_{n,k-1} - u_{n,k-2}}{\tau_n}    \right\|^2_{L^2}
= \tfrac12 \| v_{n,k} - v_{n,k-1}   \|^2_{L^2} .
$$
Since $\mathcal{J}_{n,k}$ is Fr\'echet differentiable, by minimality of $u_{n,k}$ in \eqref{e.scheme} the following Euler-Lagrange equation is easily obtained:
\begin{align}
	\tau_n^{-1} \rho \, \langle v_{n,k} - v_{n,k-1} , \phi  \rangle + \eta \langle \nabla v_{n,k} , \nabla \phi \rangle + \mu\langle \nabla u_{n,k} , \nabla \phi \rangle - ( f_{n,k} , \phi )_\U + \partial_u \Psi  ( u_{n,k} , \xi_{n,k-1} ) [ \phi ] & =  0 , \label{e.ED1}  
\end{align}
for all $\phi\in \U$ and $k=1,...,n$.
Moreover, by Remark \ref{r.sharp} we have also
\begin{align}
	\tau_n^{-1} \rho \, \langle v_{n,k} - v_{n,k-1} , \phi  \rangle  + \eta \langle \nabla v_{n,k} , \nabla \phi \rangle + \mu\langle \nabla u_{n,k} , \nabla \phi \rangle - ( f_{n,k} , \phi )_\U + \partial_u \Psi  ( u_{n,k} , \xi_{n,k} ) [ \phi ] & =  0 , \label{e.ED1bis}
\end{align}
for every $\phi \in \U$, and $k=1,...,n$. 
Finally, for every $k=1,...,n$ we have
\begin{equation} \label{e.ED3}
| \jumpo{u_{n,k}} | \le \xi_{n,k} , 
\quad
(\xi_{n,k}-\xi_{n,k-1})( | \jumpo{u_{n,k}} | -\xi_{n,k})=0, 
\quad 
| \xi_{n,k} - \xi_{n,k-1} |  \le | \jumpo{u_{n,k}}  - \jumpo{u_{n,k-1}} | ,
\end{equation} 
$\mathcal H^1$-a.e.~on $K$.

We introduce $u_n$ and $v_n$ as the piece-wise affine interpolant of $u_{n,k}$ and $v_{n,k}$, respectively, in the points $t_{n,k}$, for $k \ge 0$.  Namely
\begin{align*}
 &u_n(t):=u_{n,k-1}+(t-t_{n,k-1}) \tau_n^{-1} (u_{n,k}-u_{n,k-1}), \quad \text{ for }t\in[t_{n,k-1},t_{n,k}),
\\
 &v_n(t):=v_{n,k-1}+(t-t_{n,k-1})\tau_n^{-1}(v_{n,k}-v_{n,k-1}), \quad \text{ for }t\in[t_{n,k-1},t_{n,k}),
\\
  &\xi_n(t):=\xi_{n,k-1}+(t-t_{n,k-1}) \tau_n^{-1} (\xi_{n,k}-\xi_{n,k-1}), \quad \text{ for }t\in[t_{n,k-1},t_{n,k}).
\end{align*}
Clearly, $\dot{u}_{n,k} (t) = v_{n,k}$ for $t \in (t_{n,k} , t_{n,k+1})$.
Let us also introduce the piece-wise constant left-continuous interpolant $u^\sharp_{n}$, $v^\sharp_{n}$ and $\xi_n^\sharp$, i.e.
\begin{align}
 &u^\sharp_n(t):=u_{n,k}  
 \qquad v^\sharp_n(t):=v_{n,k} 
 \qquad \xi^\sharp_n (t):=\xi_{n,k} \quad \text{ for }t\in (t_{n,k-1},t_{n,k}]. 
\end{align}
We are ready to prove the first a-priori estimates on $u_n$ and $\xi_n$.

\begin{proposition} \label{p.compy} Let $\bar \xi>0$ and $(u_n,\xi_n)$ be given by \eqref{e.scheme}. \ERIC Under the assumptions of Theorem \ref{t.teor} the sequence $(u_n, \xi_n)$ is bounded in $W^{1,2}(0,T; \U \times \Xi)$ while the sequence $v_n$ is bounded in $W^{1,2} (0,T; \U^*)$. More precisely, there is a constant $C>0$, independent of $\tau_n$ and $\bar\xi$, such that 
\begin{align} \label{est.eps0n}
 \|u_n\|_{W^{1,2}(0,T; \U)}+\|\xi_n\|_{W^{1,2}(0,T; \Xi)}+\|v_n\|_{W^{1,2} (0,T; \U^*)}\leq C.
\end{align}
Moreover for any $1\leq p< \infty$ there is a constant $C>0$, independent of $\tau_n$ and $\bar\xi$, such that 
\begin{align} \label{est.eps1n}
&\|u_n\|_{W^{1,\infty}(0,T;L^2(\Omega))}+\|\xi_n\|_{W^{1,2} (0,T; L^p (K))}+\|v_n\|_{L^\infty(0,T; L^2(\Omega))}\leq C.
\end{align}
\end{proposition}



\proof We denote by $C$ a generic positive constant, independent of $n$ and $\bar\xi$, which might change from line to line.
By \eqref{e.ED1bis} with  $\phi=u_{n,k} - u_{n,k-1} = \tau_n v_{n,k} $  we get, for $k \ge 1$, 
\begin{align*}
	\rho \langle v_{n,k} - v_{n,k-1} , v_{n,k} \rangle + \tau_n \eta \| \nabla v_{n,k} \|^2_{L^2} & + \mu \langle \nabla u_{n,k} , \nabla (u_{n,k} - u_{n,k-1}) \rangle \, = \\[3pt] &  = \tau_n \, ( f_{n,k} , v_{n,k} )_\U- \tau_n \, \partial_u \Psi ( u_{n,k} , \xi_{n,k}) [v_{n,k}].
\end{align*}
Recalling  \eqref{e.derctrl} we have $| \partial_w \psi ( \jumpo{u} , \xi)[\jumpo{\phi}]| \le \hat\psi' ( 0) \, | \jumpo{\phi} |$, hence by continuity of the traces we get 
\begin{equation} \label{e.rem} | \partial_u \Psi ( u_{n,k} , \xi_{n,k})[ v_{n,k} ]| \le \int_K | \partial_w \psi (\jumpo{u_{n,k}} , \xi_{n,k})[ \jumpo{v_{n,k}}] | \, d\mathcal{H}^{1}  \le  \hat\psi'(0) \, \| \jumpo{v_{n,k}} \|_{L^1} \le C \| \nabla v_{n,k} \|_{L^2} .
\end{equation}
Hence
\begin{align*}
	\rho \langle v_{n,k} - v_{n,k-1} , v_{n,k} \rangle + \tau_n \eta \| \nabla v_{n,k} \|^2_{L^2} & + \mu \langle \nabla u_{n,k} , \nabla (u_{n,k} - u_{n,k-1}) \rangle - \tau_n C \| \nabla v_{n,k} \|_{L^2}   \le \tau_n \, ( f_{n,k} , v_{n,k}  )_\U .
\end{align*}
Using the scalar inequality $ \tfrac12 \eta d^2 - C_\eta \le \eta d^2 - C d$ for the second and forth term and the identity $\langle a - b , a \rangle = \tfrac12 \| a \|_{L^2}^2 - \tfrac12 \| b \|_{L^2}^2 + \tfrac12 \| a-b\|_{L^2}^2$ for the first and third term in the left-hand side, we obtain
\begin{align*}
	\rho \big( \tfrac12 \| v_{n,k} \|_{L^2}^2 - \tfrac12 \| v_{n,k-1} \|_{L^2}^2  & + \tfrac12 \| v_{n,k} - v_{n,k-1} \|_{L^2}^2 \big) + \tfrac12 \tau_n \eta \| \nabla v_{n,k} \|_{L^2}^2 \, - \tau_n C_\eta +  \\[3pt] & + \mu \big( \tfrac12 \| \nabla u_{n,k} \|_{L^2}^2 -  \tfrac12 \| \nabla u_{n,k-1} \|_{L^2}^2+   \tfrac12 \| \nabla (u_{n,k} - u_{n,k-1} ) \|_{L^2}^2  \big) \\[3pt] & \le \tau_n \,( f_{n,k} , v_{n,k} )_\U . 
\end{align*}
Let us fix $m\in\{1,...,n\}$.
Taking the sum of the previous estimate for $k=1,...,m$ we get, after neglecting and re-arranging few terms,  
\begin{align}
	\tfrac12 \rho \| v_{n,m} \|^2_{L^2}  
	& + \tfrac12 \eta \int_0^{m \tau_n} \| \nabla v_n^\sharp \|^2_{L^2}\, dt + \tfrac12 \mu \| \nabla u_{n,m} \|_{L^2}^2 
\le \nonumber \\ \le & \,
 \tfrac12 \rho \| v_{n,0} \|_{L^2}^2 +   \tfrac12 \mu \| \nabla u_{n,0} \|_{L^2}^2 + \int_0^{m \tau_n} ( f^\sharp_n ,  v^\sharp_n )_\U \, dt  + C_\eta T . \label{e.aux}
\end{align}
Since $f\in W^{1,2} (0,T;\U^*)$, we write 
$$
	\int_0^{m \tau_n} ( f_n^\sharp , v_n^\sharp )_\U \, dt=\sum_{k=1}^m( f_{n,k} ,u_{n,k}-u_{n,k-1} )_\U=-\sum_{k=0}^{m-1}( f_{n,k+1}-f_{n,k} ,u_{n,k} )_\U-( f_{n,0} ,u_{0} )_\U+ ( f_{n,m} ,u_{n,m} )_\U,
$$
which is estimated from above by
\begin{align}\label{stima_2caso}
\|\dot f\|_{L^1 (\U^*)}\|u_n^\sharp\|_{L^\infty(0,T; \U)} +\|f_{n,0}\|_{\U^*}\|u_0\|_\U+\|f_{n,m}\|_{\U^*}\|u_{n,m}\|_\U\leq C \| u_n^\sharp \|_{L^\infty(0,T; \U)} . 
\end{align}
Then, going back to \eqref{e.aux}, we obtain for every $m=1,...,n$ 
$$
	 \| u_{n,m} \|_\U^2 = \| \nabla u_{n,m} \|^2_{L^2}  \le C ( 1 + \| u_n^\sharp \|_{L^\infty(0,T;\U)} )  ,
$$
and thus $u^\sharp_n$ is bounded in $L^\infty(0,T;\U)$. Since the right-hand side of \eqref{e.aux} is bounded
%
it follows that also $v_n^\sharp$ is bounded in $L^\infty(0,T; L^2(\Omega))$. Thus (for $m=n$) we obtain that $v_n^\sharp$ is bounded in $L^2(0,T; \U)$. Recalling that $\dot{u}_n=v_n^\sharp$ we conclude that $u_n$ is bounded in $W^{1,\infty}(0,T;L^2(\Omega))$ and in $W^{1,2}(0,T; \U)$. 




Let us go back  to \eqref{e.ED1bis} and write it as 
$$
 \eta \langle \nabla v_{n,k} , \nabla \phi \rangle + \mu \langle \nabla u_{n,k} , \nabla \phi \rangle - ( f_{n,k} , \phi ) + \partial_u \Psi ( u_{n,k} , \xi_{n,k} ) [  \phi ]  =  \langle \tau_n^{-1} \rho ( v_{n,k} - v_{n,k-1} ) , - \phi  \rangle  .
$$
As $(v_{n,k} - v_{n,k-1}) \in L^2 (\Omega)$ the duality between $\U$ and $\U^*$ is represented by the $L^2$ scalar product, thus taking the supremum with respect to $\| \phi \|_{\U} \le 1$ and using again \eqref{e.derctrl} we get %
$$
	\rho \| \tau_n^{-1} ( v_{n,k} - v_{n,k-1} ) \|_{\U^*} \le \eta \| \nabla v_{n,k} \|_{L^2} + \mu \| \nabla u_{n,k} \|_{L^2} + C \| f_{n,k} \|_{\U^*} + C .
$$
Hence, by the boundedness of $u_n$ and $f_n$, we get 
$$
	\int_0^T \| \dot{v}_n \|^2_{\U^*} \, dt \le C' \int_0^T \| \nabla \dot{u}_n \|^2_{L^2} +  \| \nabla u_n \|^2_{L^2}  + \| f_n \|^2_{\U^*} + 1 \, dt < \bar{C} ,
$$
and thus $\dot{v}_n$ is bounded in $L^2 ( 0,T ; \U^*)$.

It remains to prove the compactness of the internal variable $\xi_n$. By \eqref{e.ED3} we have $| \xi_{n,k} - \xi_{n,k-1} |  \le | \jumpo{u_{n,k}}  - \jumpo{u_{n,k-1}} | $ a.e.~in $K$. Then, by the embedding of $H^{1/2} (K)$ in $L^p(K)$ (for any $1 \le p< \infty$) and the continuity of traces we have 
$$
	\| \xi_{n,k} - \xi_{n,k-1} \|_{L^p} \le \| \jumpo{u_{n,k}} - \jumpo{u_{n,k-1}} \|_{L^p} \le C \| u_{n,k} - u_{n,k-1} \|_{\U} .
$$
It follows that $\xi_n$ is bounded in $W^{1,2} (0, T; L^p(K))$. \qed 

\begin{remark} \label{r.cost1} \normalfont By \eqref{e.aux}, the initial data enter into the estimates  \eqref{est.eps0n} and \eqref{est.eps1n} only by means of the $L^2$-norms of $v_0$ e $\nabla u_0$.
\end{remark}

\subsection{Improved compactness}

We are now ready to refine the estimates provided by the last Proposition, in view of the proof of Theorem \ref{d.wsol2}. To this aim we introduce the notation  $w_{n,k} := (v_{n,k} - v_{n,k-1}) / \tau_n = \dot{v}_{n,k}$.

\begin{proposition} \label{p.kkomp} Under the assumptions of Theorem \ref{d.wsol2} there exists $C>0$ (independent of $\bar\xi$ and $n$) and $\bar{\tau} >0 $ (depending on $\bar\xi$) such that for $\tau_n < \bar{\tau}$ we have
\begin{align}\label{est.kkomp}
\| v_{n} \|_{L^\infty(0,T; \U)} + \| \dot{v}_{n} \|_{L^\infty(0,T;L^2(\Omega))}  \le C . 
\end{align}
\end{proposition}

\proof {\bf Step I.} 
 Choosing $\phi = v_{n,k} - v_{n,k-1}=\tau_n w_{n,k}$ the Euler-Lagrange equation \eqref{e.ED1bis} provides, for every $k=1,...,n$,  
\begin{align*}
  \tau_n \rho \, \langle w_{n,k}, w_{n,k} \rangle & + \eta \langle \nabla v_{n,k} , \nabla v_{n,k} - \nabla v_{n,k-1} \rangle + \mu\langle \nabla u_{n,k} , \nabla v_{n,k} - \nabla v_{n,k-1} \rangle \,+ \\ & + \partial_u \Psi  ( u_{n,k} , \xi_{n,k} ) [ v_{n,k} - v_{n,k-1}] =  \langle f_{n,k} , v_{n,k} - v_{n,k-1} \rangle,
\end{align*}
and for every $k=2,...,n$ 
\begin{align*}
  \tau_n\rho \, \langle w_{n,k-1} , w_{n,k} \rangle & + \eta \langle \nabla v_{n,k-1} , \nabla v_{n,k} - \nabla v_{n,k-1} \rangle  + \mu\langle \nabla u_{n,k-1} , \nabla v_{n,k} - \nabla v_{n,k-1} \rangle \, + \\ &+  \partial_u \Psi  ( u_{n,k-1} , \xi_{n,k-1} ) [ v_{n,k} - v_{n,k-1}] = \langle f_{n,k-1} , v_{n,k} - v_{n,k-1} \rangle .
\end{align*}
{We have $u_{n,-1} = u_{0} - \tau_n v_{0}$ and we set} $v_{n,-1} = v_{n,0}-\tau_n w_0$, where $w_0$ is the function appearing in \eqref{e.w0}. Then, for $k=1$ 
the above identity reads
\begin{align*}
\tau_n\rho\langle w_0,w_{n,1}\rangle &+ \eta \langle \nabla v_{n,0} , \nabla v_{n,1} - \nabla v_{n,0} \rangle  +\mu\langle \nabla u_{n,0} , \nabla v_{n,1}-\nabla v_{n,0}  \rangle \\
&+  \partial_u \Psi  ( u_{n,0} , \xi_{n,0} ) [ v_{n,1}-v_{n,0}] = \langle f_{n,0} , v_{n,1} -v_{n,0}\rangle ,
\end{align*}
which holds true thanks to hypothesis \eqref{e.w0}.
Hence, taking the difference of the two  identities above gives, for every $k=1,...,n$, 
\begin{align}
  \tau_n \rho \, \langle w_{n,k} - w_{n,k-1} , w_{n,k} \rangle & + \eta \| \nabla v_{n,k} - \nabla v_{n,k-1} \|^2_{L^2(\Omega)} + \tau_n \mu \langle \nabla v_{n,k} , \nabla v_{n,k} - \nabla v_{n,k-1} \rangle \,+ \nonumber  \\ & +  \partial_u \Psi  ( u_{n,k} , \xi_{n,k} ) [ v_{n,k} - v_{n,k-1}] -   \partial_u \Psi  ( u_{n,k-1} , \xi_{n,k-1} ) [ v_{n,k} - v_{n,k-1}] \nonumber \\ &   = \tau_n \langle \dot{f}_{n,k} , v_{n,k} - v_{n,k-1} \rangle  . \, \label{e.ELdif}
\end{align}
In the sequel we will employ the visco-elastic term to control the interface integral (by Poincar\'e and trace inequalities) writing
\begin{equation}  \label{e.viscotrace}
     \eta \| \nabla v_{n,k} - \nabla v_{n,k-1} \|^2_{L^2(\Omega)} \ge  \int_K c_\eta \jumpo{ v_{n,k} - v_{n,k-1}}^2 \, d\mathcal{H}^{1}.  
\end{equation} 
Now we argue a.e.~on $K$ in order to write the cohesive term in a more convenient way. Remember that $0 < \bar\xi \le \xi_{n,k}$. Denote $c_{n,k} = c_{\xi_{n,k}} = \hat\psi'(\xi_{n,k}) / \xi_{n,k}$, then by \eqref{e.Psiceps} we can write 
%
%
\begin{align*} 
	 \partial_u \Psi  ( u_{n,k} , \xi_{n,k} )  [ v_{n,k} - v_{n,k-1}] & - \partial_u \Psi  ( u_{n,k-1} , \xi_{n,k-1} ) [ v_{n,k} - v_{n,k-1}] =  \nonumber \\ 
	& = \int_K   ( c_{n,k} \jumpo{ u_{n,k} } - c_{n,k-1} \jumpo{ u_{n,k-1} }   ) \, \jumpo{v_{n,k} - v_{n,k-1}}  \, d\mathcal{H}^{1} .
\end{align*}
In order to estimate the above integral it is useful to factor out a velocity term, writing
$$
\left( c_{n,k} \jumpo{ u_{n,k} } - c_{n,k-1} \jumpo{ u_{n,k-1} }  \right) = \alpha_{n,k} \jumpo{ u_{n,k} - u_{n,k-1} } = \tau_n \alpha_{n,k} \jumpo{ v_{n,k} } ,
$$
where $\alpha_{n,k}$ is defined in the following way: if $\jumpo{u_{n,k}} = \jumpo{u_{n,k-1}}$ then $\alpha_{n,k} = c_{n,k} = c_{\xi_{n,k}}$; if $\jumpo{u_{n,k}} \neq \jumpo{u_{n,k-1}}$ then 
$$
\alpha_{n,k} = \frac{c_{n,k} \jumpo{ u_{n,k} } - c_{n,k-1} \jumpo{ u_{n,k-1}}}{\jumpo{ u_{n,k} - u_{n,k-1} }}.
$$
In particular, for $k=0$ we have $\alpha_{n,0} = c_{n,0} = c_{\xi_{n,0}}$, because $\jumpo{v_0}=0$ and thus $\jumpo{u_{n,0}} = \jumpo{u_{n,-1}}$. In this way we can always write 
\begin{align} \label{star}
	 \partial_u \Psi  ( u_{n,k} , \xi_{n,k} )  [ v_{n,k} - v_{n,k-1}] & - \partial_u \Psi  ( u_{n,k-1} , \xi_{n,k-1} ) [ v_{n,k} - v_{n,k-1}] =
\nonumber \\ & = \int_K  \tau_n \alpha_{n,k} \jumpo{ v_{n,k} } \jumpo{v_{n,k} -v_{n,k-1}} \,  d \mathcal{H}^1,
\end{align}
for every $k=1,...,n$.
Let $\bar{c} = c_{\bar\xi}>0$, $\beta>0$ be the constant appearing in \eqref{stima_beta}, and 
\begin{align}\label{alphabar_def}
 m =\max \{ k :  \xi_{n,k} \le \xi_c \} , \qquad \bar{\alpha}_{n,k} = \begin{cases}  \alpha_{n,k} & \text{if $k < m$,} \\ - \beta & \text{if $k \geq m$.} \end{cases} 
\end{align}
We define $\bar\tau := c_\eta / ( \beta + \bar{c})$ and $\gamma :=  (   \hat\psi' (0) + \xi_c \beta )^2 / \xi_c^2 c_\eta$. We claim that for $\tau_n < \bar\tau$  and for every index $k = 1, ..., n$ it holds, a.e.~on $K$, 
\begin{align}  \label{e.rud}
	  c_\eta \jumpo{v_{n,k} -v_{n,k-1}}^2 & + \tau_n \alpha_{n,k} \jumpo{ v_{n,k} } \jumpo{v_{n,k} -v_{n,k-1}}  \ge \nonumber \\ & \ge   \, \tfrac12 \tau_n  \bar{\alpha}_{n,k} \, \jumpo{v_{n,k}}^2 - \tfrac12 \tau_n  \bar{\alpha}_{n,k-1}  \jumpo{v_{n,k-1}}^2   - \tau_n^2  \gamma \, \jumpo{v_{n,k}}^2 .
\end{align}
The proof of  the above estimate is contained in the next step. 
 
\medskip
{\bf Step II.}  First, we prove that 
\begin{equation}\label{stima_alphabeta}
 - \beta \le \alpha_{n,k} \le \bar{c} = c_{\bar{\xi}} .
\end{equation}
Note that $c_{n,k} \le c_{n,k-1}$ as $\xi_{n,k} \ge \xi_{n,k-1} >0$ and $\hat\psi'(\xi_{n,k}) \le \hat\psi'(\xi_{n,k-1})$. In the case $\xi_{n,k} = \xi_{n,k-1}$ we have  $c_{n,k} = c_{n,k-1} = \alpha_{n,k} $, hence $-\beta < 0 \le c_{n,k} \le c_{n,0} \le \bar{c} $ by the monotonicity of $c_{n,k}$. In the case $\xi_{n,k} > \xi_{n,k-1}$ we have $| \jumpo{u_{n,k}}|=\xi_{n,k}>\xi_{n,k-1} \ge \jumpo{u_{n,k-1}}$. Assume that $ \jumpo{u_{n,k}}=\xi_{n,k}>0$ (the case $\jumpo{u_{n,k}}=-\xi_{n,k}$ is similar),
hence $\jumpo{u_{n,k}}-\jumpo{u_{n,k-1}}>0$. The upper bound in \eqref{stima_alphabeta} is easily achieved because $c_{n,k} \jumpo{u_{n,k}} \le c_{n,k-1} \jumpo{u_{n,k}}$, hence 
\begin{align*}
 \alpha_{n,k} = \frac{c_{n,k} \jumpo{u_{n,k}} -c_{n,k-1} \jumpo{u_{n,k-1}}}{ \jumpo{u_{n,k}}-\jumpo{u_{n,k-1}}} \le c_{n,k-1} \le \bar{c} .
\end{align*}
Let us discuss the lower bound in \eqref{stima_alphabeta}. Suppose first $\xi_{n,k}\leq \xi_c$. As $\jumpo{u_{n,k}}=\xi_{n,k}$ we have $c_{n,k} \jumpo{u_{n,k}} = \hat\psi'(\xi_{n,k})$ while $c_{n,k-1} \jumpo{u_{n,k-1}} \le c_{n,k-1} \xi_{n,k-1} = \hat \psi'(\xi_{n,k-1})$. Moreover  $\hat{\psi}' ( \xi_{n,k}) - \hat{\psi}' ( \xi_{n,k-1} )  \le 0 $ and $\jumpo{u_{n,k}}-\jumpo{u_{n,k-1}} \ge \xi_{n,k} - \xi_{n,k-1}$. Therefore by \eqref{stima_beta} we have
\begin{align}
 \alpha_{n,k} = \frac{c_{n,k} \jumpo{u_{n,k}}-c_{n,k-1} \jumpo{u_{n,k-1}}}{ \jumpo{u_{n,k}}-\jumpo{u_{n,k-1}}}& 
\geq\frac{\hat\psi'(\xi_{n,k})-\hat\psi'(\xi_{n,k-1})}{ \jumpo{u_{n,k}}-\jumpo{u_{n,k-1}}}
 \geq\frac{\hat\psi'(\xi_{n,k})-\hat\psi'(\xi_{n,k-1})}{\xi_{n,k}-\xi_{n,k-1}}\geq-\beta.\nonumber
\end{align}
Suppose now $\xi_{n,k}\geq \xi_c$ and $\xi_{n,k-1}<\xi_c$. We have $\jumpo{u_{n,k}} = \xi_{n,k} \ge \xi_c $ and $\jumpo{u_{n,k-1}} \le \xi_{n,k-1}$ and thus $ \jumpo{u_{n,k}}-\jumpo{u_{n,k-1}} \ge \xi_c - \xi_{n,k-1}$. Note that $\hat\psi' ( \xi_{n,k} ) = \hat\psi' ( \xi_c) = 0$. 
As in the previous case we write
\begin{align}
 \alpha_{n,k} = \frac{c_{n,k} \jumpo{u_{n,k}}-c_{n,k-1} \jumpo{u_{n,k-1}}}{\jumpo{u_{n,k}}-\jumpo{u_{n,k-1}}}&\geq\frac{\hat\psi'(\xi_{n,k})-\hat\psi'(\xi_{n,k-1})}{\jumpo{u_{n,k}}-\jumpo{u_{n,k-1}}}\geq \frac{\hat\psi'(\xi_c)-\hat\psi'(\xi_{n,k-1})}{\xi_c-\xi_{n,k-1}}\geq-\beta.\nonumber
\end{align}
If $\xi_{n,k-1}\geq\xi_c$ then trivially $\alpha_{n,k}=0$. Hence \eqref{stima_alphabeta} is achieved.

\medskip
Next, we prove that  for $\tau_n <\bar \tau$ and for $k \le m$ (see \eqref{alphabar_def}) it holds 
\begin{equation} \label{e.pre-rot}
	  c_\eta \jumpo{v_{n,k} -v_{n,k-1}}^2 + \tau_n \alpha_{n,k} \jumpo{ v_{n,k} } \jumpo{v_{n,k} -v_{n,k-1}}  \ge \tfrac12 \tau_n   \alpha_{n,k} \, \jumpo{v_{n,k}}^2 - \tfrac12 \tau_n  \alpha_{n,k-1} \,   \jumpo{v_{n,k-1}}^2 .
\end{equation}
Let us start writing
\begin{align*}
	  \tau_n \alpha_{n,k} \jumpo{ v_{n,k} } \jumpo{v_{n,k} -v_{n,k-1}}  & = \tfrac12 \tau_n   \alpha_{n,k} \, \jumpo{v_{n,k}}^2 - \tfrac12 \tau_n  \alpha_{n,k} \,   \jumpo{v_{n,k-1}}^2 + \tfrac12 \tau_n  \alpha_{n,k} \, \jumpo{v_{n,k} - v_{n,k-1} }^2  \\
	& \ge  \tfrac12 \tau_n   \alpha_{n,k} \, \jumpo{v_{n,k}}^2 - \tfrac12 \tau_n  \alpha_{n,k-1} \,   \jumpo{v_{n,k-1}}^2 + \tfrac12 \tau_n   ( \alpha_{n,k-1} -\alpha_{n,k} ) \jumpo{v_{n,k-1}}^2 + \\	
& \quad  - \tfrac12 \tau_n   \beta \, \jumpo{v_{n,k} - v_{n,k-1} }^2   .
\end{align*}
For $\tau_n < c_\eta/ \beta$, the last term is clearly controlled by $\tfrac12 c_\eta \jumpo{v_{n,k} - v_{n,k-1} }^2  $. We will prove that for $\tau_n < c_\eta / (\beta + \bar{c})$ we have 
\begin{align}\label{subclaim}
 \tfrac12 c_\eta \jumpo{v_{n,k} - v_{n,k-1} }^2+\tfrac12 \tau_n   ( \alpha_{n,k-1} -\alpha_{n,k} ) \jumpo{v_{n,k-1}}^2 \geq0, 
\end{align}
which implies \eqref{e.pre-rot}. 
It is not restrictive to assume that $\jumpo{v_{n,k-1}} \neq 0$, otherwise there is nothing to prove. By \eqref{stima_alphabeta} we have $\alpha_{n,k-1} - \alpha_{n,k} \ge -  \beta - \bar{c}$. We distinguish two cases:
$$
	\jumpo{v_{n,k}} \, \jumpo{v_{n,k-1}} \le 0 
	\ \mbox{ and } \ 
	\jumpo{v_{n,k}} \, \jumpo{v_{n,k-1}} > 0 .
$$
In the first case,  for $\tau_n < c_\eta/ (\beta + \bar{c})$ we have 
$$
       c_\eta  \jumpo{ v_{n,k}  - v_{n,k-1} }^2 + \tau_n   ( \alpha_{n,k-1} -\alpha_{n,k} ) \jumpo{v_{n,k-1}}^2 \ \ge  c_\eta  \, \jumpo{ v_{n,k-1} }^2  -   \tau_n  ( \beta + \bar{c}) \jumpo{v_{n,k-1}}^2 \ge 0 .
$$
The proof of the second case is more delicate. We will show that $(\alpha_{n,k-1} - \alpha_{n,k} ) \ge 0$, which implies \eqref{subclaim} and concludes the proof of \eqref{e.pre-rot}.  Assume that $\jumpo{v_{n,k-1}} > 0$, so  that $\jumpo{v_{n,k}}  > 0 $ (the case $\jumpo{v_{n,k-1}} < 0$ is treated similarly);
under this assumption
\begin{equation*} 
- \xi_{n,k-2} \le \jumpo{u_{n,k-2}} < \jumpo{u_{n,k-1}} <  \jumpo{u_{n,k}} \le \xi_{n,k} \le \xi_c . 
\end{equation*}
Note that strict inequalities hold because speeds are strictly positive, while the latter inequality holds because $k \le m$ (see \eqref{alphabar_def}). 
We need to further distinguish four cases.

\begin{itemize}

\item[1.] If $\xi_{n,k-2}=\xi_{n,k-1} = \xi_{n,k}$ then we have trivially $c_{n,k-2}=c_{n,k-1}=c_{n,k}$ and thus $\alpha_{n,k-1}=\alpha_{n,k}$.

\item[2.] If $\xi_{n,k-2} < \xi_{n,k-1} < \xi_{n,k}$ then 
$\jumpo{u_{n,k-1}} = \xi_{n,k-1}$ and $\jumpo{u_{n,k}} = \xi_{n,k}$, whereas $ \jumpo{u_{n,k-2}} \leq \xi_{n,k-2}$. Then  $c_{n,k-1} \jumpo{u_{n,k-1}} = \hat\psi'(\xi_{n,k-1})$ and $c_{n,k-2} \jumpo{u_{n,k-2}} \le \hat\psi'(\xi_{n,k-2})$, thus
\begin{align*}
      \alpha_{n,k-1} &= \frac{ c_{n,k-1} \jumpo{u_{n,k-1}}  -  c_{n,k-2} \jumpo{u_{n,k-2}} }{\jumpo{u_{n,k-1} - u_{n,k-2}}} 
      \geq \frac{ \hat\psi'(\xi_{n,k-1})-\hat\psi'(\xi_{n,k-2})}{\xi_{n,k-1}-\jumpo{u_{n,k-2}}} \\ & \geq\frac{ \hat\psi'(\xi_{n,k-1})-\hat\psi'(\xi_{n,k-2})}{\xi_{n,k-1}-\xi_{n,k-2}},
\end{align*}
while
$$
  \alpha_{n,k} = \frac{ \hat\psi'(\xi_{n,k})-\hat\psi'(\xi_{n,k-1})  }{\xi_{n,k}-\xi_{n,k-1}},
$$
and the inequality $\alpha_{n,k-1} - \alpha_{n,k} \ge 0$ follows by the concavity of $\hat\psi'$ on $[0,\xi_c]$.
%

\item[3.] If $\xi_{n,k-2} = \xi_{n,k-1} < \xi_{n,k}$ then $- \xi_{n,k-1} \le \jumpo{u_{n,k-2}} < \jumpo{u_{n,k-1}} \le \xi_{n,k-1}$ while $\jumpo{u_{n,k}} = \xi_{n,k}$.  Then $c_{n,k-1} = c_{n,k-2}$ and we can write 
$$
      \alpha_{n,k-1} = \frac{ c_{n,k-1} \jumpo{u_{n,k-1}}  -  c_{n,k-2} \jumpo{u_{n,k-2}} }{\jumpo{u_{n,k-1} - u_{n,k-2}}} =c_{n,k-1} .
$$
Since $\jumpo{u_{n,k}} = \xi_{n,k}>0$ by the monotonicity of $c_{n,k}$ we get 
$$
  \alpha_{n,k} = \frac{ c_{n,k} \jumpo{u_{n,k}}  -  c_{n,k-1} \jumpo{u_{n,k-1}}  }{\jumpo{u_{n,k} - u_{n,k-1}}}\leq \frac{ c_{n,k-1} \jumpo{u_{n,k}}  -  c_{n,k-1} \jumpo{u_{n,k-1}}  }{\jumpo{u_{n,k} - u_{n,k-1}}}=c_{n,k-1} .
$$
We conclude  $\alpha_{n,k-1} \ge \alpha_{n,k}$.


\item[4.]  The case $\xi_{n,k-2}<\xi_{n,k-1} = \xi_{n,k}$ actually does not occur because $\jumpo{u_{n,k}} \le \xi_{n,k} = \xi_{n,k-1} = \jumpo{u_{n,k-1}}$ and thus $\jumpo{v_{n,k}} \le 0$.
\end{itemize}
The proof of \eqref{subclaim} is achieved, and  \eqref{e.pre-rot} follows.

Let us now analyse the case $k \ge m+1$. We will prove that 
\begin{align} \label{e.pre-rot-bis1}
	  c_\eta \jumpo{v_{n,k} -v_{n,k-1}}^2 + \tau_n \alpha_{n,k} \jumpo{ v_{n,k} } \jumpo{v_{n,k} -v_{n,k-1}}  \ge &  - \, \tfrac12 \tau_n   \beta \, \jumpo{v_{n,k}}^2 + \tfrac12 \tau_n  \beta \,   \jumpo{v_{n,k-1}}^2   - \tau_n^2  \gamma \, \jumpo{v_{n,k}}^2 .
\end{align}

Before proceeding note that $c_{n,k} \jumpo{u_{n,k}} = 0$ for every $k > m$. Note also that for $k=m+1$ and $k=m+2$ estimate \eqref{e.pre-rot} could fail, since it may happen that $\alpha_{n,k}  =0$ and $\alpha_{n,k-1} < 0$ (this is related to the loss of concavity of $\hat\psi'$ in $\xi_c$). On the contrary, for $k>m+2$ we have $\alpha_{n,k-1}=\alpha_{n,k}=0$.

We first show that for $k > m$ we have $\alpha_{n,k} \le  \hat{\psi}'(0) / \xi_c $. Indeed, assume that $\jumpo{u_{n,k}} > \xi_c$ (the case $\jumpo{u_{n,k}} < - \xi_c$ is treated in a similar way) and
note that $c_{n,k} \jumpo{u_{n,k}} = 0$. If $\jumpo{u_{n,k-1}} \ge \xi_c$ then trivially $c_{n,k-1} \jumpo{u_{n,k-1}} = 0$ and $\alpha_{n,k}=0$.  If $0 \le \jumpo{u_{n,k-1}} \le \xi_c$ then $c_{n,k-1} \jumpo{u_{n,k-1}} \ge 0$ and thus 
$$
      \alpha_{n,k} = \frac{ c_{n,k} \jumpo{u_{n,k}}  - c_{n,k-1} \jumpo{u_{n,k-1}} }{\jumpo{ u_{n,k}  - u_{n,k-1} } } \le 0
$$
 (remember that $\jumpo{ v_{n,k} } > 0$).
If $\jumpo{u_{n,k-1}} \le 0$ then $c_{n,k-1} \jumpo{u_{n,k-1}} \le 0$ and  $ - c_{n,k-1} \jumpo{u_{n,k-1}} \le \hat\psi' ( \xi_{n,k-1}) \le \hat\psi'(0)$; since in this case  $\jumpo{ u_{n,k}  - u_{n,k-1} } \ge \xi_c$ we conclude that 
\begin{equation}\label{e.laz}
      \alpha_{n,k} = \frac{ c_{n,k} \jumpo{u_{n,k}}  - c_{n,k-1} \jumpo{u_{n,k-1}} }{\jumpo{ u_{n,k}  - u_{n,k-1} } } \le
\frac{ \hat{\psi}'(0) }{\xi_c} .  
\end{equation}
We now turn back to \eqref{e.pre-rot-bis1} and  estimate separately the terms 
\begin{gather*}
\tfrac 12 c_\eta \jumpo{v_{n,k} -v_{n,k-1}}^2 + \tau_n ( \alpha_{n,k} + \beta)  \jumpo{ v_{n,k} } \jumpo{v_{n,k} -v_{n,k-1}}  , 
\\
\tfrac12 c_\eta \jumpo{v_{n,k} -v_{n,k-1}}^2 - \tau_n \beta \jumpo{ v_{n,k} } \jumpo{v_{n,k} -v_{n,k-1}}  .
\end{gather*}
For the first we use the algebraic inequality $ \tfrac12 c_\eta s^2 +  b s \ge  \tfrac14 c_\eta s^2 - c$, with $c =b^2 / c_\eta $, which provides
\begin{align*}
\tfrac 12 c_\eta \jumpo{v_{n,k} -v_{n,k-1}}^2 + \tau_n ( \alpha_{n,k} + \beta)  \jumpo{ v_{n,k} } \jumpo{v_{n,k} -v_{n,k-1}} &  \ge \tfrac 14 c_\eta \jumpo{v_{n,k} -v_{n,k-1}}^2 - \tau^2_n \gamma \jumpo{ v_{n,k} }^2 ,
\end{align*}
since by \eqref{e.laz} 
$$
   \frac{(  \alpha_{n,k} + \beta )^2}{c_\eta} \le \frac{ (\hat\psi(0) + \xi_c \beta)^2}{\xi^2_c c_\eta}  = \gamma .
$$
For the second term, we write instead
\begin{align*}
	  - \tau_n \beta \jumpo{ v_{n,k} } \jumpo{v_{n,k} -v_{n,k-1}}  & = - \tfrac12 \tau_n   \beta \, \jumpo{v_{n,k}}^2 + \tfrac12 \tau_n  \beta \,   \jumpo{v_{n,k-1}}^2 - \tfrac12 \tau_n  \beta \, \jumpo{v_{n,k} - v_{n,k-1} }^2 .
\end{align*}
Once again the last term is balanced by $\tfrac12 c_\eta \jumpo{v_{n,k} - v_{n,k-1} }^2$ for $\tau_n < c_\eta / \beta$. 
In conclusion \eqref{e.pre-rot-bis1} follows.

%
%
%
%
%
%
%
%
%
%
%
We now conclude the proof of the claim \eqref{e.rud}. Recalling \eqref{alphabar_def}, inequality \eqref{e.pre-rot} implies that, for $k < m$, we have 
\begin{align*} 
	  c_\eta \jumpo{v_{n,k} -v_{n,k-1}}^2 & + \tau_n \alpha_{n,k} \jumpo{ v_{n,k} } \jumpo{v_{n,k} -v_{n,k-1}}  \ge \nonumber \\ & \ge   \, \tfrac12 \tau_n  \bar{\alpha}_{n,k} \, \jumpo{v_{n,k}}^2 - \tfrac12 \tau_n  \bar{\alpha}_{n,k-1}  \jumpo{v_{n,k-1}}^2   - \tau_n^2  \gamma \, \jumpo{v_{n,k}}^2 .
\end{align*}
The same inequality holds for $k \ge m+1$, thanks to \eqref{e.pre-rot-bis1}. It remains to consider 
the case $k=m$. Since $\alpha_{n,m} \ge -\beta = \bar{\alpha}_{n,m}$ and  $\alpha_{n,m-1} = \bar{\alpha}_{n,m-1}$, by \eqref{e.pre-rot} we have 
\begin{align*} 
	  c_\eta \jumpo{v_{n,m} -v_{n,m-1}}^2 & + \tau_n \alpha_{n,m} \jumpo{ v_{n,m} } \jumpo{v_{n,m} -v_{n,m-1}}  \ge \nonumber \\ & \ge   \, \tfrac12 \tau_n  \alpha_{n,m} \, \jumpo{v_{n,m}}^2 - \tfrac12 \tau_n \alpha_{n,m-1}  \jumpo{v_{n,m-1}}^2   \\ & \ge   \, - \tfrac12 \tau_n  \beta \, \jumpo{v_{n,m}}^2 - \tfrac12 \tau_n  \alpha_{n,m-1}  \jumpo{v_{n,m-1}}^2   - \tau_n^2  \gamma \, \jumpo{v_{n,m}}^2 . \\
	   & { =} \, \tfrac12 \tau_n  \bar{\alpha}_{n,m} \, \jumpo{v_{n,m}}^2 - \tfrac12 \tau_n  \bar{\alpha}_{n,m-1}  \jumpo{v_{n,m-1}}^2   - \tau_n^2  \gamma \, \jumpo{v_{n,m}}^2 .
\end{align*}

\medskip
{\bf Step III.} We now go back to \eqref{e.ELdif} and estimate the other terms.
For the inertial one we write
\begin{align}
	 \, \tau_n \rho \langle w_{n,k} - w_{n,k-1} , w_{n,k} \rangle & = \tfrac12  \tau_n\rho \|  w_{n,k} \|^2_{L^2(\Omega)} - \tfrac12 \tau_n \rho  \|  w_{n,k-1} \|^2_{L^2(\Omega)} +   \tfrac12 \tau_n \rho \| w_{n,k} - w_{n,k-1} \|^2_{L^2(\Omega)} 
\nonumber \\
	& =  \tfrac12  \tau_n \rho \|  \dot{v}_{n,k} \|^2_{L^2(\Omega)} - \tfrac12 \tau_n \rho  \|  \dot{v}_{n,k-1} \|^2_{L^2(\Omega)} +   \tfrac12 \tau_n \rho \| \dot{v}_{n,k} - \dot{v}_{n,k-1} \|^2_{L^2(\Omega)},  \label{e.1ELdif} 
	\end{align}
whereas	
	\begin{align}
	 \, \tau_n \mu \langle \nabla v_{n,k} , \nabla v_{n,k} - \nabla v_{n,k-1} \rangle & = \tfrac12  \tau_n \mu \| v_{n,k} \|^2_{\U} - \tfrac12 \tau_n \mu \| v_{n,k-1} \|^2_{\U} +   \tfrac12 \tau_n \mu \| v_{n,k} - v_{n,k-1} \|^2_{\U}. \label{e.2ELdif} 
\end{align}
The power term reads 
\begin{align}
	 \tau_n \langle \dot{f}_{n,k} , v_{n,k} - v_{n,k-1} \rangle  =  \tau_n^2  \langle \dot{f}_{n,k} , \dot{v}_{n,k} \rangle .
	\label{e.3ELdif}
\end{align}
%
For convenience, let us introduce the notation
$$
	[ v_{n,k} ]_K =  \int_K  \bar{\alpha}_{n,k} \, \jumpo{v_{n,k}}^2 \, d\mathcal{H}^{1} 
$$
(we recall that  $- \beta \le \bar\alpha_{n,k} \le \bar{c}$).
By  \eqref{e.viscotrace}, \eqref{star}, and  \eqref{e.rud} we can write 
\begin{align}
	\eta \| \nabla v_{n,k} - \nabla v_{n,k-1} \|^2_{L^2(\Omega)} & + 
		\partial_u \Psi  ( u_{n,k} , \xi_{n,k} )  [ v_{n,k} - v_{n,k-1}]  - \partial_u \Psi  ( u_{n,k-1} , \xi_{n,k-1} ) [ v_{n,k} - v_{n,k-1}]  \nonumber \\ 
   & \geq c_\eta \int_K \jumpo{ v_{n,k} - v_{n,k-1}}^2 +  ( c_{n,k} \jumpo{ u_{n,k} } - c_{n,k-1} \jumpo{ u_{n,k-1} }   ) \, \jumpo{v_{n,k} - v_{n,k-1}}  \, d\mathcal{H}^{1}  \nonumber \\ &  \ge \tfrac12 \tau_n
\int_K  \bar{\alpha}_{n,k} \, \jumpo{v_{n,k}}^2 -   \bar{\alpha}_{n,k-1}  \jumpo{v_{n,k-1}}^2 \, d\mathcal{H}^{1}   - \tau_n^2 \gamma \int_K  \jumpo{v_{n,k}}^2 \, d\mathcal{H}^{1} \nonumber \\
& =   \tfrac12 \tau_n   [ v_{n,k} ]_K -   \tfrac12 \tau_n   [ v_{n,k-1} ]_K - \tau_n^2 \gamma \| \jumpo{v_{n,k}} \|^2_{L^2(K)} .\label{e.pre-rot-bis}
\end{align}
Recalling \eqref{e.unicoe} and that $-\beta \le \bar\alpha_{n,k}$, we have
$$
     c \| \jumpo{v_{n,k}} \|^2_{L^2(K)}  \le  \tfrac12 \mu \| v_{n,k} \|^2_{\U} - \beta \| \jumpo{ v_{n,k} } \|^2_{L^2}   \le   \mu \| v_{n,k} \|^2_{\U} + [ v_{n,k} ]_K .
$$
Therefore, using (\ref{e.1ELdif}-\ref{e.pre-rot-bis}),  equation \eqref{e.ELdif} (divided by $\tau_n$) yields, for $k=1,\dots,n$,
\begin{align*} 
& \tfrac12 \big(  \rho \|  \dot{v}_{n,k} \|^2_{L^2(\Omega)} +   \mu \| v_{n,k} \|^2_{H^1} +  [ v_{n,k} ]_K  \big)
- \tfrac12 \big(   \rho  \|  \dot{v}_{n,k-1} \|^2_{L^2(\Omega)}  + \mu \| v_{n,k-1} \|^2_{H^1} + [ v_{n,k-1} ]_K \big)
\\
&\le 
 \tau_n  ( \dot{f}_{n,k} , \dot{v}_{n,k} )_{L^2(\Omega)} +    \tau_n  \gamma \, \| \jumpo{v_{n,k}} \|^2_{L^2(K)}
\\
&\le
 \tfrac12  \tau_n  \frac{1}{\rho\gamma'} \| \dot{f}_{n,k} \|^2_{L^2} +   \tau_n  \gamma' \tfrac12 (  \rho \| \dot{v}_{n,k} \|^2_{L^2} + \mu \| v_{n,k} \|^2_{H^1} + [ v_{n,k} ]_K  ),
\end{align*} 
where $\gamma'= 2\gamma / c$. 
Remember that  $ v_{n,-1}=v_{n,0}-\tau_nw_0$, hence $\dot{v}_{n,0}=w_0$ and $\jumpo{v_{n,0}}=0$; by Gronwall inequality we conclude that there exists $\gamma''>0$ such that for every $1 \le k \le n$ and for $\tau_n \leq \bar\tau = c_\eta / (\beta + \bar{c})$ it holds
\begin{align*}
   \rho \|  \dot{v}_{n,k} \|^2_{L^2(\Omega)} +   \mu \| v_{n,k} \|^2_{H^1} +  [ v_{n,k} ]_K & \le \gamma'' ( \rho \|  \dot{v}_{n,0} \|^2_{L^2(\Omega)} +   \mu \| v_{n,0} \|^2_{H^1} +  [ v_{n,0} ]_K )  +  \gamma'' \| f \|^2_{W^{1,2}(0,T; L^2)}  \le C_0.
\end{align*}
Note that $\gamma''$ depends only on $\gamma =  (   \hat\psi' (0) + \xi_c \beta )^2 / \xi_c^2 c_\eta$, in particular it is independent of $\bar\xi$. Note also that the upper bound $C_0$ is independent of $\bar\xi$ since $[ v_{n,0} ]_K = 0$.
Invoking again  \eqref{e.unicoe} the thesis follows. \qed

\begin{remark} \label{r.cost2} \normalfont As far as the initial data, the constant $C$ appearing in \eqref{est.kkomp} depends only on the $L^2$-norms of $w_0 $ e $\nabla v_0$. 
\end{remark}

\subsection{Existence}

In this section we pass to the limit as $\tau_n\rightarrow 0$, or equivalently as $n\rightarrow+\infty$, to obtain an evolution $(u ,\xi)$ which solves the system \eqref{maineq} and which enjoys good compactness properties. 

\begin{theorem} \label{p.teor2} Let $\bar{\xi} > 0$ and $\xi_0 \ge \bar\xi$ a.e.~on $K$.  

Under the assumptions of Theorem \ref{t.teor}  there exists an evolution $(u ,\xi) \in W^{1,2} (0,T; \U \times \Xi)$ with $\dot{u} \in W^{1,2} (0,T;\U^*)$ such that for a.e.~$t \in (0,T)$
\begin{align}\label{maineq2}
	\begin{cases}
	\rho \ddot u  (t) + \partial_u \F ( t , u  (t) , \xi  (t) ) + \partial_v \mathcal{R} ( \dot{u} (t) ) = 0 ,  & \text{in $\U^*$,} \\[5pt]
	\dot{\xi}  (t)  ( \xi  (t) -  | \jumpo{u  (t)} | ) = 0, \ \dot{\xi}  (t) \ge 0 , \text{ and }  \ | \jumpo{u (t)} | \le \xi  (t) ,   & \text{a.e.~in $K$,} \\[4pt]
	u  (0) = u_0 , \ \dot{u}  (0) = v_0, \ { \xi(0)=\xi_0 }.
	\end{cases}
\end{align}
Moreover   for any $1\leq p< \infty$ there exists a constant $C>0$, independent of $\bar\xi$, such that 
\begin{align}
 &\|u \|_{W^{1,2}(0,T; \U)}+\|\xi \|_{W^{1,2}(0,T; \Xi)} + \|u \|_{W^{2,2}(0,T;\U^*)} \leq C,\label{est.eps00}\\
&\|u \|_{W^{1,\infty}(0,T;L^2(\Omega))}+ \|\xi \|_{W^{1,2} (0,T; L^p (K))}\leq C. \label{est.eps11}
\end{align}
Under the stronger assumptions of Theorem \ref{d.wsol2} there exists also $C>0$, independent of $\bar\xi$, such that 
 \begin{align}\label{est.eps22}
 \| \dot u \|_{L^\infty(0,T;\U)} + \| \ddot{u} \|_{L^\infty(0,T; L^2(\Omega))}  \le C.  
 \end{align}
\end{theorem}



\proof {\bf Step I.} Let us prove \eqref{est.eps00}-\eqref{est.eps22}. By Proposition \ref{p.compy} we know that (up to subsequences) $( u_n , \xi_n ) \weakto ( u , \xi)$ in $W^{1,2}(0,T; \U \times \Xi)$ and $v_n \weakto v$ in $W^{1,2} (0,T; \U^*)$. 
It follows that $u_n (t_n) \weakto u(t)$ in $\U$ if $t_n \to t$ and thus $u^\sharp_n (t) \weakto u(t)$ in $\U$, because $u_n^\sharp (t) = u_n (t_{n,k_n+1})$ for $t \in ( t_{n,k_n} , t_{n,k_n+1}]$. Since $\dot{u}_n (t) = v_n^\sharp(t)$, Lemma \ref{a.l.int} and the embedding of $\U$ in $\U^*$ yield 
$$	\| \dot{u}_n (t) - v_n (t) \|_{\U^*} = \| v_n^\sharp (t) - v_n (t) \|_{\U^*} \le C \tau_n^{1/2}  ,
\qquad \dot{u}_n \weakto \dot{u} \quad \text{ in $L^2(0,T, \U^*)$.} $$
Since $v_n \weakto v$ in $L^2(0,T, \U^*)$, it follows that $\dot{u} = v$ in $\U^*$. Then, inequality \eqref{est.eps00} follows from \eqref{est.eps0n} by weak lower semi-continuity of the norms.
Upon extracting a further subsequence (non-relabelled) we get also inequality \eqref{est.eps11} from \eqref{est.eps1n} and inequality \eqref{est.eps22} from \eqref{est.kkomp}.


Being $u_n(0) = u_0$, $\xi(0)=\xi_0$, and $v_n(0)=v_0$, it is easy to check that the initial conditions in \eqref{maineq2} are satisfied. 


\medskip
\noindent {\bf Step II.} Let us first check the Karush-Kuhn-Tucker conditions for $\xi$, partially adapting the proof of \cite[Theorem 5.8]{NegriVitali_IFB18}. Since $\xi_n \weakto \xi$ in $W^{1,2}(0,T;L^2(K))$ we have $\xi_n (t) \weakto \xi(t)$ in $L^2(K)$ for every $t \in [0,T]$. We also have $u_n (t) \weakto u (t)$ in $H^1 (\Omega)$, thus by the compact embedding of the traces we get $\jumpo{u_n(t)} \to \jumpo{u(t)}$ in $L^2(K)$ for all $t \in [0,T]$. Moreover, $u_n$ is bounded in $L^\infty (0,T;H^1(\Omega))$ and thus $\jumpo{u_n}$ is bounded in $L^\infty(0,T; L^2(K))$. Thus, by dominated convergence, we obtain that $\jumpo{u_n} \to \jumpo{u}$ in $L^2(0,T;L^2(K))$. Moreover, since $\jumpo{u_n}$ is bounded in $W^{1,2}(0,T;L^2(K))$ invoking Lemma \ref{a.l.int} we have that  $\jumpo{u_n^\sharp} \to \jumpo{u}$ in $L^2(0,T;L^2(K))$. 

From \eqref{e.ED3} for a.e.~$t\in (0,T)$ we have 
$$
    \dot{\xi}_n (t) \ge 0 , \qquad | \jumpo{u_{n}^\sharp  (t) } | \le \xi_{n}^\sharp (t) , 
\qquad
\dot{\xi}_{n} (t) ( | \jumpo{u_{n}^\sharp(t)} | -\xi_{n}^\sharp(t)) = 0 .
$$
As $\xi_n$ is monotone non-decreasing it turns out, by weak convergence, that $\xi$ is monotone non-decreasing. Since $u_n$ and $\xi_n$ are affine interpolant, from $| \jumpo{u_{n}^\sharp  (t) } | \le \xi_{n}^\sharp (t)$ we easily get $| \jumpo{u_{n}  (t) } | \le \xi_{n} (t)$. Hence, by weak convergence, $| \jumpo{u (t) } | \le \xi (t)$. As a consequence $\dot{\xi} (t) ( | \jumpo{u (t) } | - \xi (t) ) \le 0$ for a.e.~$t\in (0,T)$ and a.e.~in $K$; thus, to conclude the proof of the Karush-Kuhn-Tucker condition, it will be enough to show that 
\begin{align}\label{claim}
	\int_0^T \langle \dot{\xi} (t) , | \jumpo{u (t) } | - \xi (t)  \rangle \, dt \ge 0 .
\end{align}
We aim to pass to the limit the identity
$$
	 \int_0^T \langle  \dot{\xi}_{n} (t) ,  | \jumpo{u_{n}^\sharp(t)} | \rangle  \, dt =  \int_0^T \langle \dot{\xi}_{n} (t) , \xi_{n}^\sharp(t) \rangle \, dt .
$$
Since, $\jumpo{u_n^\sharp} \to \jumpo{u}$ (strongly) in $L^2(0,T; L^2(K))$ and  $\dot{\xi}_n \rightharpoonup \dot{\xi}$ in $L^2(0,T; L^2(K))$ we get 
\begin{equation}\label{e.l1}
\int_0^T \langle \dot{\xi}_n(t) , |\jumpo{u_n^\sharp(t)}| \rangle \, dt  \rightarrow \int_0^T \langle \dot{\xi}(t) , |\jumpo{u(t)}| \rangle \, dt ; 
\end{equation}
in order to conclude \eqref{claim}, it suffices to prove that 
$$ \liminf_{n \to +\infty} \int_0^T \langle \dot{\xi}_n(t) , \xi_n^\sharp (t) \rangle \, dt  \ge \liminf_{n \to +\infty} \int_0^T \langle \dot{\xi}_n(t) , \xi_n (t) \rangle \, dt  \ge  \int_0^T \langle \dot{\xi}(t) , \xi(t) \rangle \, dt , $$ 
where, in the first inequality, we have used the fact that $\dot{\xi}_n \ge 0$ and $\xi_n^\sharp \ge \xi_n $. By Lions-Magenes lemma (see e.g. \cite[Lemma 1.2 Ch. III \S 1]{Temam-NS:1977} 
) we have 
$$
	\int_0^T \langle \dot{\xi}_n(t) , \xi_n (t) \rangle \, dt  = 
\tfrac12 \| \xi_n (T) \|^2_{L^2} - \tfrac12 \| \xi_0 \|^2_{L^2} .
$$
Since $\xi_n(T) \weakto \xi (T)$ in $L^2(K)$ we have 
\begin{align}
	\liminf_{n \to +\infty} \int_0^T \langle \dot{\xi}_n(t) , \xi_n (t) \rangle \, dt  & \ge \liminf_{n \to +\infty} \, \tfrac12 \| \xi_n (T) \|^2_{L^2} - \tfrac12 \| \xi_0 \|^2_{L^2} \nonumber \\ & \ge \| \xi (T) \|^2_{L^2} - \tfrac12 \| \xi_0 \|^2_{L^2} = \int_0^T \langle \dot{\xi} (t) , \xi (t) \rangle \, dt .\label{e.l2}
\end{align}
The claim and the Karush-Kuhn-Tucker conditions are proved. 

\medskip
\noindent {\bf Step III.} We prove that $\xi_n \rightarrow \xi$ strongly in $L^2(0,T;L^2(K))$ (using the argument of \cite[Lemma 4.14]{NScala_NARWA17}): since $\dot{\xi} \, \xi = \dot{\xi} \, | \jumpo{u} |$ and 
$\dot{\xi}_n \, \xi^\sharp_n = \dot{\xi}_n \, | \jumpo{u_n^\sharp} |$ a.e.~in $(0,T)$, using \eqref{e.l1} and arguing as in \eqref{e.l2} for every time $t \in (0,T)$ we can write
\begin{align*} 
\tfrac{1}{2}\|\xi(t)\|_{L^2 }^2&=\tfrac{1}{2}\|\xi_0 \|_{L^2 }^2+\int_0^t\langle \dot{\xi}(r),| \jumpo{u(r)}|\rangle \,dr  \\ 
& =\tfrac{1}{2}\|\xi_0\|_{L^2 }^2+\lim_{n\rightarrow\infty}\int_0^t\langle \dot{\xi}_n(r),| \jumpo{u_n^\sharp (r)} |\rangle \, dr \\
& = \tfrac{1}{2}\|\xi_0\|_{L^2 }^2+\lim_{n\rightarrow\infty}\int_0^t\langle \dot{\xi}_n(r), \xi_{n}^\sharp (r)\rangle \, dr \\ 
&   \ge \tfrac{1}{2}\|\xi_0\|_{L^2 }^2+ \limsup_{n\rightarrow\infty} \int_0^t\langle \dot{\xi}_n(r), \xi_n (r)\rangle \, dr  \geq\limsup_{n\rightarrow\infty}\tfrac{1}{2}\|\xi_n(t)\|_{L^2 }^2 .
 \end{align*}
On the other hand, $\xi_n (t) \rightharpoonup \xi(t)$ in $L^2(K)$ for every $t \in (0,T)$, thus $\| \xi(t)  \|^2_{L^2} \le  \liminf_{n \to +\infty} \| \xi_n (t) \|^2_{L^2}$. It follows that $\| \xi_n (t) \|_{L^2} \to \| \xi(t)  \|_{L^2}$ and thus $\xi_n (t) \to \xi(t)$ in $L^2(K)$. We conclude by dominated convergence because $\xi_n$ is bounded in $L^\infty(0,T;L^2(K))$. It follows that $\xi^\sharp_n \rightarrow \xi$ strongly in $L^2(0,T;L^2(K))$.

\medskip
\noindent {\bf Step IV.} From \eqref{e.ED1bis} for every $0 \le t_1 < t_2 < T$ we get 
\begin{equation} \label{e.part}
	\int_{t_1}^{t_2} (\rho \dot{v}_n (t) , \phi )_\U + \langle \mu\nabla u^\sharp_n (t) + \eta\nabla \dot{u}_n (t) , \nabla \phi \rangle - ( f^\sharp_n (t) , \phi )_\U + \partial_u \Psi ( u_n^\sharp (t) , \xi^\sharp_n (t)) [\phi]   \, dt =0,  
\end{equation}
for every variation $\phi \in \U$. As $\dot{v}_n \weakto \dot{v} = \ddot{u}$ in $L^2(0,T; \U^*)$ we have 
$$
	\int_{t_1}^{t_2} ( { \rho}\dot{v}_n (t) , \phi )_\U \, dt \to \int_{t_1}^{t_2} ( {\rho}\dot{v} (t) , \phi )_\U \, dt = \int_{t_1}^{t_2} ( {\rho}\ddot{u} (t) , \phi )_\U \, dt \,.
$$
Since $u_n^\sharp$ is bounded in $L^\infty (0,T; \U)$ (cf.~Proposition \ref{p.compy}) and $u^\sharp_n (t) \weakto u(t)$ in $\U$, we infer 
$$
	\int_{t_1}^{t_2} \langle \mu\nabla u^\sharp_n (t) , \nabla \phi \rangle \, dt \to \int_{t_1}^{t_2}  \langle\mu \nabla u (t), \nabla \phi \rangle \, dt \,.
$$
Since $u_n \weakto u$ in $W^{1,2} (0,T; \U)$ and $f^\sharp_n \to f$ in $L^2(0,T;\U^*)$ we get immediately that 
$$
	\int_{t_1}^{t_2} \langle\eta \nabla \dot{u}_n (t) , \nabla \phi \rangle - ( f^\sharp_n (t) , \phi )_\U  \, dt \to \int_{t_1}^{t_2}  \langle\eta \nabla \dot{u} (t), \nabla \phi \rangle - (f (t) , \phi )_\U\, dt \,.
$$
Let us see that 
\begin{align}\label{claim_psieps}
	\lim_{n  \to +\infty} \int_{t_1}^{t_2} \partial_u \Psi ( u^\sharp_n (t) , \xi^\sharp_n (t) ; \phi) \, dt =  \int_{t_1}^{t_2}\partial_u \Psi ( u (t) , \xi (t) ; \phi) \, dt. 
\end{align}
From the previous steps, we know that $\jumpo{u^\sharp_n} \to u$ and $\xi^\sharp_n \to \xi$ in $L^2 (0,T ; L^2(K))$. Thus, we extract a subsequence (non-relabelled) such that the convergence of $\jumpo{u^\sharp_n}$ and $\xi^\sharp_n$ holds pointwise $\mathcal{H}^1$-a.e. on $K$ and a.e. on $(0,T)$. Recalling that $\xi_n\geq \xi_0\geq \bar \xi>0$, we infer that a.e.~in $K\times[0,T]$ it holds
\begin{align}
 \partial_u \psi ( \jumpo{u^\sharp_n}, \xi^\sharp_n ; \jumpo{\phi}) =  \partial_u \psi ( \jumpo{u^\sharp_n }, \xi^\sharp_n ) \jumpo{\phi}\rightarrow \partial_u \psi ( \jumpo{u} , \xi) \jumpo{\phi}=  \partial_u \psi ( \jumpo{u}, \xi; \jumpo{\phi}),
\end{align}
where we have used the continuity of $\partial_u\psi$ away from $(0,0)$. Finally, by the fact that $\jumpo{u_n}\leq \xi_n$, we also get the bound
$$
|\partial_u \psi ( \jumpo{u_n} , \xi_n )| \leq \tilde \psi'(\xi_n)\leq \tilde \psi'(0),
$$
so that \eqref{claim_psieps} follows by dominated convergence.

In conclusion, taking the limit in \eqref{e.part} we get 
$$
	\int_{t_1}^{t_2} (\rho \ddot{u}(t) , \phi )_\U + \langle \mu\nabla u (t) + \eta\nabla \dot{u} (t) , \nabla \phi \rangle -( f (t) , \phi )_\U+ \partial_u \Psi ( u(t) , \xi (t) ; \phi )  \, dt = 0, 
$$
and by arbitrariness of $t_1$ and $t_2$ we obtain the equilibrium equation. 

\qed



\section{Proofs of the main results}\label{section5}

We fix a small parameter $\bar\xi_\eps>0$  with $\bar\xi_\eps \to 0^+$ for $\eps \to 0^+$.
We truncate (from below) the initial condition $\xi_0$ with $\bar\xi_\eps$ in order to gain regularity on the energy and then on the solution $u^\eps$. In a second step we will pass to the limit as $\eps \to 0^+$ to obtain a solution to the original Cauchy problem.


Invoking Theorem \ref{p.teor2} we get the following result.

\begin{corollary} \label{p.cor2}  Assume (H1), (H2), let $(u_0 , \xi_0) \in \U \times \Xi$ with $| \jumpo{u_0} | \le \xi_0$, let $v_0\in L^2(\Omega)$, and let $f \in W^{1,2} (0,T;\U^*)$. Then for every $\eps>0$ there exists an evolution $(u^\eps,\xi^\eps) \in W^{1,2} (0,T; \U \times \Xi)$ with $\dot{u}^\eps \in W^{1,2} (0,T;\U^*)$ such that for a.e.~$t \in (0,T)$
\begin{align}\label{maineq2bis}
	\begin{cases}
	\rho \ddot u^\eps (t) + \partial_u \F ( t , u^\eps (t) , \xi^\eps (t) ) + \partial_v \mathcal{R} ( \dot{u}^\eps(t) ) = 0 ,  & \text{in $\U^*$,} \\[5pt]
	\dot{\xi}^\eps (t)  ( \xi^\eps (t) -  | \jumpo{u^\eps (t)} | ) = 0, \ \dot{\xi}^\eps (t) \ge 0 , \text{ and }  \ | \jumpo{u^\eps(t)} | \le \xi^\eps (t) ,   & \text{a.e.~in $K$,} \\[4pt]
	u^\eps (0) = u_0 , \ \dot{u}^\eps (0) =v_0 ,\  \xi^\eps (0) = \widehat\xi^\eps_0=\max\{ \bar\xi_\eps , \xi_0 \} .
	\end{cases}
\end{align}
%
Moreover, for every $1\leq p< \infty$ there exists a constant $C>0$, independent of $\eps>0$, such that 
\begin{align}
 &\|u^\eps\|_{W^{1,2}(0,T; \U)}+\|\xi^\eps\|_{W^{1,2}(0,T; \Xi)} + \|u^\eps\|_{W^{2,2}(0,T;\U^*)}\leq C,\nonumber\\
&\|u^\eps\|_{W^{1,\infty}(0,T;L^2(\Omega))}+\|\xi^\eps\|_{W^{1,2} (0,T; L^p (K))}
\leq C. \label{est.eps1}
\end{align}
 
\end{corollary}

\subsection{Solutions with higher time regularity}

Assume now  the additional hypotheses (H3) and (H4), and that $f\in W^{1,2}(0,T;L^2(\Om))$. We want to show the counterpart of Corollary \ref{p.cor2} for solutions with higher regularity. Assume that  the initial data $u_0,v_0\in \U$, $\xi_0\in \Xi$, satisfy $\jumpo{v_0}=0$ on $K$, $|\jumpo{u_0}|\leq \xi_0$ on $K$. First of all, note that, even if the initial data $(u_0, v_0, \xi_0)$ enjoy the equilibrium condition \eqref{e.w0}, i.e.,
$$
\rho w_0  + \partial_u \F ( t , u_0 (t) , \xi_0 (t) ) + \partial_v \mathcal{R} ( {v}_0)  \ni0 ,  
$$
 in general the same condition does not hold for the regularized initial data $u_0,v_0,\widehat \xi_0^\eps$, therefore, in order to apply Proposition \ref{p.kkomp}
we have first to modify $u_0,v_0,\widehat \xi_0^\eps$ in a suitable way.
To this aim, let $w_0\in L^2(\Omega)$ be the function appearing in condition \eqref{e.w0}. We define
\begin{align}\label{def_initialdataeps}
&v_0^\eps:=v_0,\nonumber\\
&u_0^\eps\in \argmin \{\mathcal F(0, u,\widehat\xi^\eps_0)+\langle\eta\nabla v^\eps_0,\nabla u\rangle+\langle \rho w_0,u\rangle, \;u\in \U \},\\
&\xi^\eps_0:= {\max\{\widehat\xi_0^\eps,|\jumpo{u^\eps_0}|\}}. \nonumber
\end{align}
Notice that the above minimization problem does not provide an initial datum $u_0^\eps$ which satisfies, in general, the constraint $|\jumpo{u_0^\eps}|\leq \widehat \xi_0^\eps$. To fix this issue we have updated the internal variable a posteriori. 
We claim that the equilibrium condition \eqref{e.w0}
 holds for the triple $(u_0^\eps,v_0^\eps, \xi_0^\eps)$. To show this we will see that $u_0^\eps$ is also a solution of the minimum problem
\begin{align}\label{minimalityofinitialdatum}
 u_0^\eps\in \argmin \{\mathcal F(0, u,\xi^\eps_0)+\langle\eta\nabla v^\eps_0,\nabla u\rangle+\langle \rho w_0,u\rangle, \;u\in \U \}.
\end{align}
To prove \eqref{minimalityofinitialdatum}, as in Remark \ref{r.sharp}, it is sufficient to observe that   thanks to the properties of the cohesive potential, one has
\begin{align}
 \psi(\jumpo{u_0^\eps}, \xi_0^\eps)= {\psi (\jumpo{u_0^\eps},\widehat\xi_0^\eps),}
 \qquad\text{ and }\qquad\psi(\jumpo{u}, \widehat\xi_0^\eps)\leq \psi(\jumpo{u},\xi_0^\eps)\;\;\forall u\in \U,
\end{align}
so that 
\begin{align}
 \mathcal F(0, u_0^\eps,\xi^\eps_0)+\langle\eta\nabla v^\eps_0,\nabla u_0^\eps\rangle+\langle \rho w_0,u_0^\eps\rangle&= \mathcal F(0, u_0^\eps,\widehat \xi^\eps_0)+\langle\eta\nabla v^\eps_0,\nabla u_0^\eps\rangle+\langle \rho w_0,u_0^\eps\rangle\nonumber\\
 &\leq  \mathcal F(0, u, \widehat \xi^\eps_0)+\langle\eta\nabla v^\eps_0,\nabla u\rangle+\langle \rho w_0,u\rangle\nonumber\\
 &\leq \mathcal F(0, u, \xi^\eps_0)+\langle\eta\nabla v^\eps_0,\nabla u\rangle+\langle \rho w_0,u\rangle,
\end{align}
for all $u\in\U$,
where we have used the minimality of $u_0^\eps$ in \eqref{def_initialdataeps}. Condition \eqref{minimalityofinitialdatum} follows, and we infer
\begin{align}
 \langle \rho w_0,\phi\rangle +\partial_u\mathcal F(0,u_0^\eps,\xi^\eps_0;\phi)+\langle\eta\nabla v^\eps_0,\nabla \phi\rangle=0.
\end{align}

We have now to show that the new initial data introduced in \eqref{def_initialdataeps} suitably converge to $(u_0, v_0, \xi_0)$. For $v_0$ there is nothing to prove, for $u_0$ and $\xi_0$ we have at disposal the following result:

\begin{lemma}
 The energy functional 
 $$u\mapsto \mathcal F(0, u,\widehat \xi^\eps_0)+\langle\eta\nabla v^\eps_0,\nabla u\rangle+\langle \rho w_0,u\rangle,$$
$\Gamma$-converges, with respect to the weak topology of $\U$, to the functional
$$u\mapsto\mathcal F(0, u,\xi_0)+\langle\eta\nabla v_0,\nabla u\rangle+\langle \rho w_0,u\rangle.$$
 \end{lemma}
 The proof of this Lemma is straightforward and we drop the details; essentially, it is based on the fact that the cohesive energy term $\Psi(u,\widehat\xi_0^\eps)$ (which is the only one involving $\eps$) well behaves in the passage to the limit. Indeed, $\widehat \xi^\eps_0= \max \{ \xi_0 , \bar \xi_\eps \} \rightarrow \xi_0$ strongly in $L^2(\Omega)$ and $\jumpo{u^\eps}\rightarrow \jumpo{u}$ strongly in $L^2(K)$, by compact embedding, when $u^\eps\rightharpoonup u$ weakly in $\U$. The convergence $\Psi(u_\epsilon ,\widehat \xi_0^\eps)\rightarrow \Psi(u,\xi_0)$ follows by continuity of $\Psi$ and by dominated convergence theorem.

 Thanks to hypothesis (H4) the above functionals are strictly convex, so that by coercivity they admit a unique minimizer. From the properties of $\Gamma$-convergence we readily see that 
 \begin{align}\label{conv.init1}
  u_0^\eps\rightarrow u_0 \text{ in }\U.
 \end{align}
 Now, since again by compact embedding $| \jumpo{u^\eps} | \rightarrow | \jumpo{u}|$ strongly in $L^2(K)$, we also easily obtain
 \begin{align}\label{conv.init2}
   \xi^\eps_0= \max \{ \widehat \xi^\eps_0 , | \jumpo{ u^\eps_0} | \} \rightarrow \xi_0\text{ in }L^2(K),
 \end{align}
where we have used that $\widehat \xi_0^\eps\rightarrow \xi_0$ strongly in  $L^2(K)$.
 
 Finally, as a further consequence of the convergences above, the norms $\|u_0^\eps\|_\U$, $\|v_0^\eps\|_\U$, $\| \xi_0^\eps\|_{L^2(K)}$ are uniformly bounded as $\eps\rightarrow 0$. Therefore also the constants appearing in the right-hand side of estimates \eqref{est.eps0n}, \eqref{est.eps1n}, and \eqref{est.kkomp} are uniformly bounded (see Remark \ref{r.cost1} and \ref{r.cost2}). We hence arrive to
 
\begin{corollary}\label{p.cor2bis}
Under the hypotheses of Corollary \ref{p.cor2}, suppose in addition (H3), (H4), and that $f\in W^{1,2}(0,T;L^2(\Om))$; moreover assume that $u_0,v_0\in \U$, $\xi_0\in \Xi$, satisfy  $\jumpo{v_0}=0$ on $K$, $| \jumpo{u_0} |\leq \xi_0$ on $K$, and  condition \eqref{e.w0}. Then  for every $\eps>0$ there exists  an evolution $(u^\eps,\xi^\eps) \in W^{1,2} (0,T; \U \times \Xi)$ with $\dot{u}^\eps \in W^{1,2} (0,T;\U^*)$ satisfying
\begin{align*} 
	\begin{cases}
	\rho \ddot u^\eps (t) + \partial_u \F ( t , u^\eps (t) , \xi^\eps (t) ) + \partial_v \mathcal{R} ( \dot{u}^\eps(t) ) = 0 ,  & \text{in $\U^*$,} \\[5pt]
	\dot{\xi}^\eps (t)  ( \xi^\eps (t) -  | \jumpo{u^\eps (t)} | ) = 0, \ \dot{\xi}^\eps (t) \ge 0 , \text{ and }  \ | \jumpo{u^\eps(t)} | \le \xi^\eps (t) ,   & \text{a.e.~in $K$,} \\[4pt]
	u^\eps (0) = u^\eps_0 , \ \dot{u}^\eps (0) =v_0 ,\  \xi^\eps (0) = \xi_0^\eps , 
	\end{cases}
\end{align*}
i.e., the conclusion of Corollary \ref{p.cor2} with the initial data $u_0^\eps,v_0^\eps, \xi_0^\eps$  in \eqref{def_initialdataeps}. Moreover there is a constant $C>0$, independent of $\eps>0$, such that 
 \begin{align}\label{est.eps2}
 \| \dot u^\eps \|_{L^\infty(0,T;H^1 (\Omega))} + \| \ddot{u}^\eps \|_{L^\infty(0,T; L^2(\Omega))}  \le C.  
 \end{align}
\end{corollary}

\subsection{Proof of Theorem \ref{t.teor} and \ref{d.wsol2}}

By Lemma \ref{l.sub=var} the differential inclusion in the first equation of \eqref{maineq} in variational form reads 
\begin{align}\label{maineq3bis}
( \rho \ddot u (t) , \phi )_\U + \partial_u \F ( t , u (t) , \xi(t)  ; \phi ) + \langle \eta\nabla \dot u (t) , \nabla \phi \rangle \ge 0 ,  \ \text{for every $\phi \in \U$,}  
\end{align}
for a.e. $t\in[0,T]$.

\begin{proposition} \label{p.conv2} Along with the assumptions of Corollary \ref{p.cor2}, there exists $(u, \xi)\in W^{1,2}(0,T; \U \times \Xi) $ such that, for a non-relabelled subsequence, as $\eps\rightarrow0$ it holds 
\begin{align}
 &u^\eps\weakto u\text{ weakly in }W^{2,2}(0,T; \U^*)\cap W^{1,2}(0,T;\U)\text{ and weakly* in  }W^{1,\infty} (0,T; L^2(\Om)),\nonumber\\
 &\xi^\eps\weakto \xi\text{ weakly in }W^{1,2}(0,T; \Xi).
\label{conv.2}
\end{align}
Moreover $(u, \xi)$ is a solution to \eqref{maineq}.
\end{proposition}

\proof
The convergences in \eqref{conv.2} follow directly from the a-priori estimates \eqref{est.eps1}. 
To prove that the couple $(u,\xi)$ satisfies \eqref{maineq} we argue as in Theorem \ref{p.teor2} and show first that the Karush-Kuhn-Tucker conditions hold, namely the second line in \eqref{maineq}. Following Step II and Step III of the proof of Theorem \ref{p.teor2} we easily infer also that 
\begin{align}\label{strongxi}
 \xi^\eps\rightarrow\xi \text{ strongly in }L^{2}(0,T; L^2(K)).
\end{align}
Let us turn to the first equation in \eqref{maineq}, that is, we have to prove that for all $\phi\in \U$ equation \eqref{maineq3bis} holds for a.e. $t\in(0,T)$.
From \eqref{maineq2bis} for every $0 \le t_1 < t_2 < T$ we have 
\begin{equation}\label{e.part2}
	\int_{t_1}^{t_2} {( \rho\ddot{u}^\eps (t) , \phi )_\U} + \langle \mu\nabla u^\eps (t) + \eta\nabla \dot{u}^\eps (t) , \nabla \phi \rangle - ( f (t) , \phi )_\U + \partial_u \Psi ( u^\eps (t) , \xi^\eps (t)) [\phi]   \, dt =0. 
\end{equation}
Under \eqref{conv.2} all the terms pass to the limit but the cohesive one. For this we show that
$$
	\limsup_{\eps  \to 0} \int_{t_1}^{t_2} \partial_u \Psi ( u^\eps (t) , \xi^\eps (t) ; \phi) \, dt \le \int_{t_1}^{t_2} \partial_u \Psi ( u (t) , \xi (t) ; \phi) \, dt .
$$
Up to a  (non-relabelled) subsequence, thanks to \eqref{strongxi} we can assume that $(\jumpo{u^\eps},\xi^\eps)\rightarrow(\jumpo{u},\xi)$ pointwise almost everywhere in $(0,T)\times K$.
%
Using Lemma \ref{l.usc-psi'}, Fatou's Lemma, and Fubini's Theorem, we get 
\begin{align*}
\limsup_{\eps\to0} \int_{t_1}^{t_2}  \partial_u \Psi ( u^\eps (t) , \xi^\eps (t) ; \phi) \, dt & = \limsup_{\eps  \to 0} \int_{t_1}^{t_2}  \int_K \partial_w \psi ( \jumpo{u^\eps (t)} , \xi^\eps(t) ; \jumpo{\phi}  ) \, d\mathcal{H}^1 \, dt  \\ & \le \int_{t_1}^{t_2}  \int_K \partial_w \psi ( \jumpo{u(t)} , \xi(t) ; \jumpo{\phi}  ) \, d\mathcal{H}^1 \, dt =  \int_{t_1}^{t_2} \partial_u \Psi ( u (t) , \xi (t) ; \phi) \, dt. 
\end{align*}
In conclusion, taking the limsup in \eqref{e.part2} we get 
$$
	\int_{t_1}^{t_2} ({\rho} \ddot{u}(t) , \phi )_\U + \langle \mu\nabla u (t) + \eta\nabla \dot{u} (t) , \nabla \phi \rangle - \langle f (t) , \phi \rangle + \partial_u \Psi ( u(t) , \xi (t) ; \phi )  \, dt \ge 0 
$$
By arbitrariness of $t_1$ and $t_2$ we obtain the equilibrium inequality. 
Finally, since $u^\eps(0)=u_0$, $\xi^\eps(0)=\xi^\eps_0\rightarrow \xi_0$ in $\Xi$, and $v^\eps(0)=v_0$, from \eqref{conv.2} we easily infer that $(u,\xi)$ satisfies the initial data.
\qed
\medskip

We now aim to prove the same result for the solutions provided by Corollary \ref{p.cor2bis}. Notice that the same proof of the previous proposition applies, up to show that the modified initial data $(u_0^\eps,v_0^\eps, \xi_0^\eps)$ converge to $(u_0,v_0,\xi_0)$;
this is ensured by \eqref{conv.init1} and \eqref{conv.init2}.

\begin{proposition} \label{p.conv2bis} Under the assumptions of Corollary \ref{p.cor2bis}, there exists $(u, \xi)\in W^{1,2}(0,T; \U \times \Xi) $ such that, for a not-relabelled subsequence, \eqref{conv.2} holds and  $(u, \xi)$ is a solution to \eqref{maineq}. Moreover $$u\in W^{1,\infty}(0,T;{\U})\cap W^{2,\infty}(0,T;L^2(\Omega)).  $$
\end{proposition}

\subsection{Energy balance}

Before proceeding to the proof of the energy balance we need the following result.

\begin{lemma} \label{l.phi-phi2} Let $(u, \xi)$ be as in Proposition \ref{p.conv2}, then for a.e.~$t \in (0,T)$ it holds 
\begin{equation}  \label{e.enba}
	(\rho\ddot u (t) , \dot{u} (t) )_\U + \partial_u \F ( t , u (t) , \xi(t)  ; \dot{u} (t) ) + \langle\eta \nabla \dot u (t) , \nabla \dot{u} (t) \rangle = 0 .
\end{equation}
\end{lemma} 

\proof 
Since $u$ belongs to  $W^{1,2} (0,T; { \U})$, for a.e.~$t \in (0,T)$ we have 
$$	\lim_{h \to 0} \frac{u ( t+h) - u(t)}{h} = \dot{u}(t) \quad \text{in $H^1(\Om)$.}   $$
By continuity of the traces we have also
\begin{equation} \label{e.jumpvel}	\lim_{h \to 0}  \frac{\jumpo{u(t+h)} - \jumpo{u(t)}}{h} = \jumpo{\dot{u}(t)} \quad \text{in $L^2(K)$}.    \end{equation}
Let $t\in (0,T)$ be such that \eqref{e.jumpvel} holds. By the regularity of $u$ we have that $\dot u\in \U$ and thus
\begin{gather} \label{e.step2}
(\ddot u (t) , \pm \dot{u} (t) )_\U + \partial_u \F ( t , u (t) , \xi(t)  ; \pm \dot{u}(t) ) + ( \nabla \dot u (t) , \pm \nabla\dot{u}(t) ) \ge 0 \,.
\end{gather}

We claim that  
$$\partial_u \F ( u (t) , \xi (t) ; -\dot{u} (t)) = - \partial_u \F ( u (t) , \xi (t) ; \dot{u} (t)) , 
$$ 
which, together with \eqref{e.step2}, concludes the proof. Since the elastic energy $\E$ is differentiable it is enough to check 
$$\partial_u \Psi ( u (t) , \xi (t) ; - \dot{u}(t)) = - \partial_u \Psi ( u (t) , \xi (t) ; \dot{u}(t)). $$
We will argue pointwise in $K$. To this end, denote by $K_0 (t)$ the subset of $K$ where $\xi (t) = 0$. Consider a sequence $t_n \nearrow t$ such that 
$$	\frac{\jumpo{u(t_n)} - \jumpo{u(t)}}{t_n - t} \to  \jumpo{\dot{u}(t)} \quad \text{a.e. on $K$.}    $$
 Clearly $| \jumpo{u(t_n)} | \le \xi (t_n) \le \xi(t)$ (a.e.~on $K$ and for every $n \in \mathbb{N}$). By monotonicity of $\xi$ and by the Karush-Kuhn-Tucker conditions, a.e.~on $K_0 (t)$ we have $ | \jumpo{u(t_n)} | \le \xi (t_n) \le \xi(t) = 0$ and $| \jumpo{u(t)} | \le \xi(t) = 0$, then 
 \begin{align}\label{KKT_partial}
\dot{\xi}(t) = \jumpo{\dot{u}(t)}=0\text{ a.e.~on } K_0 (t).                                                                                                                                                    \end{align}
Remember that 
$$
	\partial_u \Psi ( u (t), \xi(t) ; \dot{u}(t) ) = \int_K \partial_w \psi \big( \jumpo{u (t,l)} , \xi(t,l) ; \jumpo {\dot{u}(t,l)}  \big) \, d\mathcal{H}^1(l) .
$$
If $l \in K \setminus K_0 (t) $ then  $\xi(t,l) > 0$; therefore the density $\psi ( \cdot , \xi (t,l))$ is differentiable in $\jumpo{u (t,l)}$ and 
$$
	\partial_w \psi ( \jumpo{u(t,l)}, \xi (t,l) ; - \jumpo{ \dot{u}(t,l)} ) = - \partial_w \psi ( \jumpo{u(t,l)}, \xi (t,l) ; \jumpo{ \dot{u}(t,l) } ) . 
$$
On the contrary, if $l \in K_0 (t)$ then $ | \jumpo{u(t,l)} | \le  \xi(t,l) = 0$ and the density $\psi ( \cdot , \xi(t,l) )$ admits only directional derivatives, however $\jumpo{ \dot{u}(t,l)} = 0$ by \eqref{KKT_partial}, and thus 
$$
	\partial_w \psi ( \jumpo{u(t,l)}, \xi (t,l) ; - \jumpo{ \dot{u}(t,l)} ) = - \partial_w \psi ( \jumpo{u(t,l)}, \xi (t,l) ; \jumpo{ \dot{u}(t,l)} ) = 0 .  
$$
The conclusion follows.
\qed

%

\begin{lemma} \label{l.PhiAC2} The map $t \mapsto \Psi ( \jumpo{u(t)} , \xi(t) )$ is absolutely continuous and 
\begin{equation} \label{e.claim}
	\Psi ( u(t^*) , \xi (t^*) ) - \Psi ( u_0 , \xi_0 ) = \int_0^{t^*} \tfrac{d}{dt} \Psi ( u(t) , \xi (t) ) \, dt = \int_0^{t^*} \partial_u \Psi \big( u(t) , \xi (t) ; \dot{u}(t) \big) \, dt,
\end{equation} 
for every $t^*\in[0,T]$.
\end{lemma}

\proof The absolute continuity follows from the fact that both $\jumpo{u}$ and $\xi$ belong to $W^{1,2} (0,T; L^2(K))$ together with the Lipschitz continuity of the cohesive energy density (see Lemma \ref{l.psi}). It follows that the map $t \mapsto \Psi ( u (t), \xi(t) )$ is absolutely continuous; thus it is a.e.~differentiable in $(0,T)$ and it is enough to show that for a.e.~$t\in(0,T)$ it holds
\begin{align}\label{claime}
	\tfrac{d}{dt} \Psi ( u(t) , \xi (t) ) = \partial_u \Psi \big( u(t) , \xi (t) ; \dot{u}(t) \big) .
\end{align}
We devide the proof of \eqref{claime} into two steps.


\noindent \textbf{Step I.} 
We claim that for a.e. $t\in(0,T)$ it holds 
\begin{equation} \label{e.KKTextra2}
   \dot{\xi} (t) \le | \jumpo{\dot{u}(t)} | \ \text{a.e.~in $K$,} \qquad 
   \dot{\xi} (t) = | \jumpo{\dot{u}(t)} |=0 \ \text{a.e.~in $K_0(t)=\{ \xi (t) = 0 \}$.} 
\end{equation}
The latter property has already been proved in Lemma \ref{l.phi-phi2}. As for the former, let us first observe that by the time regularity of $\jumpo{u}$ and $\xi$ we have 
\begin{equation} \label{e.pke}
	 \frac{\jumpo{u(t+h) - u(t)}}{h} \to \jumpo{\dot{u}(t)} \ \text{ and } \ \frac{\xi(t+h) - \xi(t)}{h} \to \dot\xi(t) 
\end{equation} 
as $h\rightarrow0$ (strongly) in $L^2 (K)$ for a.e.~$t \in (0,T)$. Moreover, for a.e.~$t \in (0,T)$ we have 
\begin{equation} \label{e.KKTaux}  \dot{\xi} (t)  ( \xi(t) -  | \jumpo{u(t)} | ) = 0 \ \text{ and } \ | \jumpo{u(t)} | \le \xi (t) \quad \text{a.e.~on $K$.} \end{equation}
We fix $t \in (0,T)$ such that both $t \mapsto \Psi ( u (t), \xi(t) )$ is differentiable and \eqref{e.pke}-\eqref{e.KKTaux} hold. 

We now prove  the first condition in \eqref{e.KKTextra2}.
Assume, by contradiction, that this does not hold and $\dot{\xi} (t) > | \jumpo{\dot{u}(t)} |$ in a set $K' \subset K$ of positive measure. Then, by \eqref{e.pke} and by convergence in measure there exists $\bar h\ll1$ such that  for every $0 < h \le \bar{h}$ we have  
$$
	\xi(t+h) - \xi(t) > | \jumpo{u(t+h) - u(t)} |  \quad \text{in $L^2(K'')$,} 
$$
for some positive measured set $K''\subset K'$.
Then, for every $0 < h \le \bar{h}$ a.e.~on $K''$ we have 
\begin{align*}
	\jumpo{u (t+h)} & \le \jumpo{u(t)} +  | \jumpo{u ( t + h )} - \jumpo{u(t)} | 
					\le \xi(t) + | \jumpo{u ( t + h )} - \jumpo{u(t)} | < \xi (t+h) .
\end{align*}
As a consequence of the Karush-Kuhn-Tucker condition $\dot{\xi} ( t+h) ( \jumpo{u(t+h)}  - \xi(t+h) ) = 0$ we have $\dot{\xi}(t+h)=0$ in $L^2(K'')$ for every $0 < h \le \bar{h}$. Thus, $\xi(t+h) = \xi(t) $ in $L^2(K'')$ and this contradicts the inequality $\xi(t+h) - \xi(t) > | \jumpo{u(t+h) - u(t)} | \ge 0$.

%
%
%

\noindent {\bf Step II.} We fix a sequence $h_n \to 0^+$  (depending on $t$) such that 
\begin{equation} \label{e.pkea.e.2}
	 \frac{\jumpo{u(t+h_n) - u(t)}}{h_n} \to \jumpo{\dot{u}(t)} \ \text{ and } \ \frac{\xi(t+h_n) - \xi(t)}{h_n} \to \dot\xi(t) \quad \text{a.e.~on $K$.}
\end{equation} 
Clearly
$$
	\tfrac{d}{dt} \Psi ( u(t) , \xi (t) ) = \lim_{n \to +\infty} \int_K \frac{\psi \big( \jumpo{u ( t + h_n ) } , \xi ( t + h_n) \big) - \psi \big( \jumpo{u ( t) } , \xi ( t) \big)}{h_n} \, d \mathcal{H}^1.
$$
Since $\psi$ is Lipschitz continuous we can write
\begin{align}
	\left| \frac{\psi \big( \jumpo{u ( t + h_n  ) } , \xi ( t + h_n ) \big) - \psi \big( \jumpo{u ( t ) } , \xi ( t ) \big)}{h_n} \right|
	\le & \ C  \left| \frac{ \jumpo{u ( t + h_n  ) }- \jumpo{u ( t ) } }{h_n} \right| \, + \nonumber \\
	& + C \left|  \frac{ \xi ( t + h_n  ) - \xi  ( t ) }{h_n} \right| \,. \label{e.lip}
\end{align}
The right-hand side converges strongly in $L^2 (K)$ and thus by dominated convergence, in order to show \eqref{claime}, it is enough to prove that 
\begin{equation} \label{e.e2}
	\lim_{n \to +\infty} \frac{\psi \big( \jumpo{u ( t + h_n ) } , \xi ( t + h_n) \big) - \psi \big( \jumpo{u ( t) } , \xi ( t) \big)}{h_n}  =  \partial_w \psi \big( \jumpo{ u (t) } , \xi (t) ;  \jumpo{ \dot{u} (t) } \big)  
\quad \text{a.e.~in $K$.}
\end{equation}

%
%
%
%
%
%
%
%
%
%

Denote $K_0(t):=\{l\in K:\xi(t,l)=0\}$ and $K_c(t):=\{l\in K:\xi(t,l) \ge \xi_c\}$. Let $l \in K$ be such that \eqref{e.KKTaux} and \eqref{e.pkea.e.2} hold. 

If $l \not \in K_0 (t) \cup K_c (t)$ then $0 < \xi ( t, l ) < \xi_c $ and thus, by Lemma \ref{l.psi}, $\psi$ is differentiable  at $( \jumpo{ u (t,l) } , \xi (t,l))$ and 
\begin{align*}
	 \lim_{n \to +\infty} \frac{\psi \big( \jumpo{u ( t + h_n , l) } , \xi ( t + h_n , l) \big) - \psi \big( \jumpo{u ( t , l) } , \xi ( t , l) \big)}{h_n} = & \ \partial_w \psi \big( \jumpo{ u (t,l) } , \xi (t,l) \big) \jumpo{ \dot{u} (t,l) } \, + \\ & +  \partial_\xi \psi \big( \jumpo{ u (t,l) } , \xi (t,l) \big) \dot{\xi} ( t, l) .
\end{align*}
By \eqref{derivative_xi} we get (for $| w | \le \xi   \neq 0$) 
$$
	\partial_\xi \psi (w, \xi   ) = - \tfrac12 \left( \frac{\hat\psi' ( \xi )}{\xi } \right)' ( \xi - w) (\xi + w) .
$$
Thus, by \eqref{e.KKTaux}
$$
	\partial_\xi \psi \big( \jumpo{ u (t,l) } , \xi (t,l) \big) \dot{\xi} ( t, l) =  - \tfrac12 \left( \frac{\hat\psi' ( \xi (t,l))}{\xi (t,l)} \right)' ( \xi (t,l) +  | \jumpo{ u (t,l) } | ) ( \xi (t,l)  -  |\jumpo{u (t,l)}|  ) \, \dot{\xi}(t,l) = 0 \,.
$$
Hence \eqref{e.e2} is proved.

%
%

If $l \in K_0(t)$ then $\xi (t, l) =0$ and by \eqref{e.KKTextra2} we have $\dot{\xi} ( t, l) = | \jumpo{ \dot{u}(t,l)} | = 0$. Remember that the cohesive density $\psi$ is no longer differentiable; however we have directional derivatives. Thus we write 
\begin{align} \label{e.differ}
	 \psi \big( \jumpo{u(t+h_n,l)} , \xi(t+h_n,l) \big) - \psi ( 0 , 0  ) 
 = & \ \psi \big( \jumpo{u(t+h_n,l)} , \xi(t+h_n,l) \big) - \psi \big( \jumpo{u(t+h_n,l)} , 0 \big) \, + \nonumber \\ 
 & +  \psi \big( \jumpo{u(t+h_n,l)} , 0 \big) - \psi \big( 0 , 0 \big) .
\end{align}
Let us first 
check that 
\begin{equation}\label{subclaime}
\lim_{n \to +\infty} \frac{\psi \big( \jumpo{u(t+h_n,l)} , \xi(t+h_n,l) \big) - \psi \big( \jumpo{u(t+h_n,l)} , 0 \big)}{h_n} = 0 . 
\end{equation}
By monotonicity (see Lemma \ref{l.psi}) 
\begin{align*}
\hat{\psi} ( \xi(t+h_n,l) ) & 
\ge \psi \big( | \jumpo{u(t+h_n,l)} |, \xi(t+h_n,l)\big) = \psi \big( \jumpo{u(t+h_n,l)} , \xi(t+h_n,l)\big)  \\ & \ge \psi \big( | \jumpo{u(t+h_n,l)} | , 0 \big) = \hat{\psi} \big( | \jumpo{u(t+h_n,l)} | \big) = \psi \big( \jumpo{u(t+h_n,l)} , 0 \big)
\end{align*}  
Hence, by the monotonicity properties of $\psi$ and the Lipschitz continuity of $\hat{\psi}$ we get 
\begin{align*}
	  \left|   \psi \big( \jumpo{u(t+h_n,l)} , \xi(t+h_n,l) \big) - \psi \big( \jumpo{u(t+h_n,l)} , 0 \big)   \right|  & =  \psi \big( \jumpo{u(t+h_n,l)} , \xi(t+h_n,l) \big) - \psi \big( \jumpo{u(t+h_n,l)} , 0 \big)  \\ & \le     \hat{\psi} ( \xi(t+h_n,l)) -  \hat{\psi} \big( | \jumpo{u(t+h_n,l)} | \big) \\ & \le c \, \big( \xi(t+h_n,l) - | \jumpo{u(t+h_n,l)} |  \big) .
\end{align*}
Moreover, being $\xi (t,l) = \jumpo{u(t,l)}=0$, we can write 
$$
	\frac{\big( \xi (t + h_n,l)  - | \jumpo{u ( t+ h_n,l)} | \big)}{h_n}  = \frac{\big( \xi (t + h_n,l) - \xi(t,l) \big) - | \jumpo{ u(t+h_n,l) - u(t,l) } |}{h_n}  \to \dot{\xi}(t,l) - | \jumpo{\dot{u} (t,l) } | = 0 \,.
$$
Thus \eqref{subclaime} is proved.
Now, let us now consider the second line in \eqref{e.differ}; it is enough to write 
$$
	\lim_{n \to +\infty} \frac{ \psi \big( \jumpo{u(t+h_n,l)} , 0 \big) - \psi \big( 0 , 0 \big) }{h_n} = \partial_w \psi \big( 0 , 0 ; \jumpo{\dot{u}(t,l)} \big) = \hat{\psi}' ( 0) | \jumpo{\dot{u}(t,l)} | = 0  .
$$
Therefore, for a.e.~$l \in K_0(t)$ it holds
$$ \lim_{n \to +\infty} \frac{\psi \big( \jumpo{u ( t + h_n , l) } , \xi ( t + h_n , l) \big) - \psi \big( \jumpo{u ( t , l) } , \xi ( t , l) \big)}{h_n} = \partial_w \psi \big( \jumpo{ u (t,l) } , \xi (t,l) ;  \jumpo{ \dot{u} (t,l) } \big). 
$$


Finally, if $l \in K_c (t)$ then $\xi ( t) \ge \xi_c$ and $\xi(t + h_n ) \ge \xi_c$ (since $h_n \ge0$) and thus 
$$
{ \psi \big( \jumpo{u ( t + h_n , l) } , \xi ( t + h_n , l) \big) = \psi \big( \jumpo{u ( t , l) } , \xi (t , l) \big) = 
\hat{\psi} ( \xi_c) .}
$$
It follows that 
$$
\lim_{n \to +\infty} \frac{\psi \big( \jumpo{u ( t + h_n , l) } , \xi ( t + h_n , l) \big) - \psi \big( \jumpo{u ( t , l) } , \xi ( t , l) \big)}{h_n} = \partial_w \psi \big( \jumpo{ u (t,l) } , \xi (t,l) ;  \jumpo{ \dot{u} (t,l) } \big) = 0.
$$
 The proof is concluded. \qed

\begin{proposition} \label{p.enbal} Let $(u, \xi)$ be as in Proposition \ref{p.conv2}, then for a.e.~$t^* \in [0,T]$ we have the energy identity
\begin{align*}
	\E ( u( t^*) ) + \Psi( u(t^*) , \xi (t^*)) + \K ( \dot{u} (t^*) )  & =  \, \E ( u_0 ) + \Psi ( u_0, \xi_0) + \K ( v_0 ) \nonumber \\ & \quad + \int_0^{t^*} ( f(t) , \dot{u} (t) )_{\U} \, dt - \int_0^{t^*} \eta \| \nabla \dot{u} (t) \|^2_{L^2} \,dt 
\end{align*}
\end{proposition}

\proof By Lions-Magenes lemma, e.g.~\cite[Lemma 1.2 Ch.~III \S 1]{Temam-NS:1977}, we can write 
$$ \tfrac12 \| \nabla{u} (t^*) \|^2_{L^2} - \tfrac12 \| \nabla{u}_0 \|^2_{L^2}  = \int_0^{t^*} \tfrac{d}{dt} ( \tfrac12 \| \nabla u (t) \|^2_{L^2}) \, dt  = \int_0^{t^*}  \langle \nabla u (t) , \nabla \dot{u} (t)  \rangle \, dt . $$
Hence, by Lemma \ref{l.PhiAC2} we have 
\begin{align*}
\E (u(t^*)) + \Psi (u(t^*) , \xi(t^*) )  & - \E (u_0) - \Psi ( u_0 , \xi_0)  = \int_0^{t^*}  \tfrac{d}{dt} \Big( \tfrac12 \mu \| \nabla u (t) \|_{L^2}^2+ \Psi ( u(t) , \xi (t) ) \Big) \, dt  \\
	& = \int_0^{t^*} \mu\langle \nabla u (t) , \nabla \dot{u}  (t)\rangle   + \partial_u \Psi \big( u(t) , \xi (t) ; \dot{u}(t) \big)  \, dt \\
	& = \int_0^{t^*} \partial_u \F ( u(t) , \xi(t) ; \dot{u}(t) )  \, dt + \int_0^{t^*} (  f(t) , \dot{u} (t) )_{\U} \, dt .
\end{align*}
Using again Lions-Magenes lemma, we get 
$$ \tfrac12 \| \dot{u} (t^*) \|^2_{L^2} - \tfrac12 \| \dot{u}_0 \|^2_{L^2}  = \int_0^{t^*} \tfrac{d}{dt} ( \tfrac12 \| \dot{u} (t) \|^2_{L^2}) \, dt   =  \int_0^{t^*}  ( \ddot{u} (t) , \dot{u}(t) )_\U \, dt     .
$$
Combining the previous two identities with 
$$
	(\rho\ddot u (t) , \dot{u} (t) )_\U + \partial_u \F ( t , u (t) , \xi(t)  ; \dot{u} (t) ) + \langle\eta \nabla \dot u (t) , \nabla \dot{u} (t) \rangle = 0 ,
$$
proved in Lemma \ref{l.phi-phi2}, we finally obtain the energy identity. \qed

\subsection{Strong solutions in polygonal domains}\label{pdes_section}

Following \cite[\S 1.5.2]{Grisvard}, we introduce the Sobolev space 
$$\widetilde{H}^{1/2} (K^+) =\{ w \in H^{1/2} ( K) :  \bar{w}^+ \in H^{1/2} ( \partial \Omega^+) \} , $$
where $\bar{w}^+$ is the null extension  of $w$ on $\partial\Omega^+$. In a similar way we introduce $\widetilde{H}^{1/2} (K^-)$. Being $K$ a polygonal domain, $\widetilde{H}^{1/2} (K^+) = \widetilde{H}^{1/2} (K^-)$. Therefore the notation $\mathcal{W} = \widetilde{H}^{1/2} (K) $ is well justified. 

We endow $\mathcal{W}$ with the standard norm of $H^{1/2}(K)$ and denote by $\mathcal{W}^* = \widetilde{H}^{-1/2} (K)$ the dual space. By \cite[Theorem 1.5.3.10]{Grisvard} whenever $\sigma \in L^2(\Omega ; \R^2)$ with $\mathrm{div}\,\stress \in L^2 (\Omega)$
we can introduce the linear and continuous Neumann operators 
\begin{align}
 & ( \stress^+ \nu , \varphi )_{\mathcal{W}} = - ( \stress^+ \nu^+ , \varphi )_{\mathcal{W}} := - \int_{\Om^+} \stress \cdot \nabla \phi^+ \, dx - \int_{\Om^+} (\mathrm{div} \, \stress ) \phi^+ \, dx,\nonumber \\ 
  & ( \stress^- \nu , \varphi )_{\mathcal{W}} = ( \stress^- \nu^- , \varphi )_{\mathcal{W}} :=\int_{\Om^-} \stress \cdot \nabla \phi^- \, dx + \int_{\Om^-} (\mathrm{div} \, \stress ) \phi^- \, dx,\nonumber  
\end{align}
where, in the left-hand side, $( \cdot,\cdot)_{\mathcal{W}}$ denotes the duality between $\mathcal{W}$ and its dual $\mathcal{W}^*$, while, in the right-hand side, $\phi^\pm \in H^1(\Omega^\pm)$ denote any liftings of the null extensions $\bar \varphi^\pm \in H^{1/2} ( \partial \Omega^\pm)$. 

In a similar way, we consider also the spaces
$ \mathcal{V}^\pm = \widetilde H^{1/2}(\partial_N\Om^\pm)$ 
and the Neumann operators 
\begin{align}
 & ( \stress \nu , \varphi )_{\mathcal V^\pm} := \int_{\Om^\pm} \stress \cdot \nabla \phi^\pm \, dx + \int_{\Om^\pm} (\mathrm{div} \, \stress ) \phi^\pm \, dx,\nonumber 
\end{align}
where 
 $\phi^\pm \in H^1(\Omega^\pm)$ denote any liftings of the null extensions $\bar \varphi^\pm \in H^{1/2} ( \partial \Omega^\pm)$.

 By Definition \ref{d.wsol} the inclusion
$$ \rho\ddot u (t) + \partial_u \F ( t , u (t) , \xi(t) ) + \partial_v \mathcal{R} ( \dot{u})  \ni0 $$ 
in variational form reads 
\begin{equation}\label{eq}
\langle \rho\ddot u(t),\phi \rangle +\langle\mu\nabla u(t)+ \eta\nabla \dot u(t),\nabla \phi\rangle+\partial \Psi(u(t),\xi(t);\phi)\geq\langle f(t),\phi\rangle,
\end{equation}
for all $\phi\in \U$.
Choosing test functions $\phi \in H^1_0 (\Omega)$ we infer that, for a.e.~$t\in[0,T]$,
\begin{equation}\label{eq0}
\langle\rho \ddot u(t),\phi \rangle +\langle  \mu\nabla u(t)+ \eta\nabla \dot u(t),\nabla \phi\rangle=\langle f(t),\phi\rangle .
\end{equation}
We easily conclude that, for a.e.~$t\in[0,T]$,
\begin{align}\label{eq1}
 \rho\ddot{u}(t) - \Delta (  \mu u(t) +  \eta\dot{u}(t) )  = f (t) \qquad& \hbox{in $H^{-1}(\Omega)$.}
\end{align}
 Since $f\in W^{1,2}(0,T;L^2(\Om))$ and since $u \in W^{2,\infty}(0,T;L^2(\Omega))$  we get that  $ \Delta ( \mu  u(t) +  \eta\dot{u}(t) )$ belongs to $L^2(0,T;L^2(\Om))$; in particular equation \eqref{eq1} holds in $L^2(\Om)$ and a.e.~on $\Om$.
Then we are allowed to define the Neumann operators $\sigma^\pm (t) \nu$ 
and $\sigma (t) \nu$ 
(see  \S\,\ref{s.poly}). 
Since $\sigma (t) = \mu\nabla u(t) + \eta\nabla \dot{u} (t)$ by \eqref{eq0} we have 
\begin{align*}
 ( \stress (t) \nu , \varphi )_{\mathcal V^+} & = \int_{\Om^+} \stress (t)  \cdot \nabla \phi^+ +  (\mathrm{div} \, \stress (t) ) \phi^+ \, dx \\ &  =  \int_{\Om^+} (\mu\nabla u(t) + \eta\nabla \dot{u} (t)) \cdot \nabla \phi^+ + (\rho\ddot{u} (t) - f(t)) \phi^+ \, dx = 0  , \nonumber 
\end{align*}
for every $\varphi \in \mathcal{V} = \widetilde{H}^{1/2} ( \partial_N \Omega^+)$. Hence $\sigma (t)  \nu =0$ in $\widetilde{H}^{-1/2} ( \partial_N \Omega^+)$.
In a similar way, for every $\varphi \in \mathcal{W} = \widetilde{H}^{1/2} ( K )$
\begin{align*}
 & ( \stress^+ \nu , \varphi )_{\mathcal{W}} = 
- \int_{\Om^+} (\mu\nabla u (t) + \eta\nabla \dot{u} (t)) \cdot \nabla \phi^+ + (\rho\ddot{u} (t) - f(t)) \phi^+ \, dx,   \\ 
  & ( \stress^- \nu , \varphi )_{\mathcal{W}} = 
  \int_{\Om^-} (\mu\nabla u(t) + \eta\nabla \dot{u} (t)) \cdot \nabla \phi^- + (\rho\ddot{u} (t) - f(t)) \phi^- \, dx . \nonumber  
\end{align*}
Therefore, for $\phi = \phi^+ 1_{\Omega^+} + \phi^- 1_{\Omega^-}$ by \eqref{eq0} we get
$$
         ( \stress^- \nu -  \sigma^+ \nu , \varphi )_{\mathcal{W}} = \int_{\Om} (\mu\nabla u(t) + \eta\nabla \dot{u} (t)) \cdot \nabla \phi + (\rho\ddot{u} (t) - f(t)) \phi \, dx = 0 ,
$$
because $\jumpo{\phi} =0$; hence $\stress^- \nu = \sigma^+ \nu$ in $\widetilde{H}^{-1/2} (K)$.


It remains to prove that $\sigma^+ (t) \nu \in \partial_w \psi ( \jumpo{u(t)}, \xi(t) )$. To this end we invoke \cite[Theorem 7.8]{NegriVitali_IFB18} where an analogous result is proved. The thesis is achieved.

\section*{Acknowledgements}
The present paper benefits of the financial support of the GNAMPA (Gruppo Nazionale per l’Analisi
Matematica, la Probabilit\`a e le loro Applicazioni) of INdAM (Istituto Nazionale di Alta Matematica).



\appendix

\section{Appendix} 

We provide here this useful classical Lemma. 

\begin{lemma} \label{a.l.int} Let $u \in W^{1,2} ( 0,T; X)$, where $X$ is a Banach space, and let $u_n^\sharp$ be its piecewise constant interpolant in the points $t_{n,k} = k T / n$. Then $\| u^\sharp_n (t) - u (t) \|_X \le \| u \|_{W^{1,2} (0,T;X)} \tau^{1/2}$ for every $t \in [0,T]$.
\end{lemma}

%
%

\bibliographystyle{plain}
\bibliography{bibliography} 

\end{document}